# Wiener-Hermite Polynomial Expansion for Multivariate Gaussian Probability Measures☆

Sharif Rahman[1]

*Applied Mathematical and Computational Sciences, The University of Iowa, Iowa City, Iowa 52242, U.S.A.*

**Abstract**

This paper introduces a new generalized polynomial chaos expansion (PCE) comprising multivariate Hermite orthogonal polynomials in dependent Gaussian random variables. The second-moment properties of Hermite polynomials reveal a weakly orthogonal system when obtained for a general Gaussian probability measure. Still, the exponential integrability of norm allows the Hermite polynomials to constitute a complete set and hence a basis in a Hilbert space. The completeness is vitally important for the convergence of the generalized PCE to the correct limit. The optimality of the generalized PCE and the approximation quality due to truncation are discussed. New analytical formulae are proposed to calculate the mean and variance of a generalized PCE approximation of a general output variable in terms of the expansion coefficients and statistical properties of Hermite polynomials. However, unlike in the classical PCE, calculating the coefficients of the generalized PCE requires solving a coupled system of linear equations. Besides, the variance formula of the generalized PCE contains additional terms due to statistical dependence among Gaussian variables. The additional terms vanish when the Gaussian variables are statistically independent, reverting the generalized PCE to the classical PCE. Numerical examples illustrate the generalized PCE approximation in estimating the statistical properties of various output variables.

## 1. Introduction

The Wiener-Hermite polynomial chaos expansion (PCE), hereafter referred to as the classical PCE, is an infinite series expansion of a square-integrable random variable involving Hermite orthogonal polynomials in independent Gaussian random variables. Introduced by Wiener [1] in conjunction with the homogeneous chaos theory, Cameron and Martin [2] proved convergence of PCE to the correct limit in the $L^2$ sense for an arbitrary random variable with finite variance.[2] Later developments include truncation of the classical PCE in a Galerkin framework, leading to a spectral stochastic finite-element method [4] and extension to a generalized PCE to account for non-Gaussian variables [5]. However, the conditions for generalization mandate completeness of measure-consistent orthogonal polynomials, as clarified only recently [3]. Approximations stemming from truncated PCE, whether classical or generalized, are commonly used for solving uncertainty quantification problems, mostly in the context of solving stochastic partial differential equations [6, 7], yielding approximate second-moment statistics of a stochastic output variable of interest. A majority of these studies, including many not cited here for brevity, address low-dimensional problems, that is, when the number of input random variables is not overly large, say, less than ten. In that case, PCE approximations can be sufficiently accurate and are known to offer significant computational advantages over crude Monte Carlo simulation (MCS), although there are exceptions [8, 9]. In high dimensions, however, PCE requires an astronomically large number of polynomials or coefficients, succumbing to the curse of dimensionality [10, 11].

---

☆Grant sponsor: U.S. National Science Foundation; Grant No. CMMI-1462385.

*Email address:* sharif-rahman@uiowa.edu (Sharif Rahman)

[1]Professor.

[2]More precisely, Cameron and Martin [2] proved convergence of PCE to the expanded random variable for a special probability space. In a more general setting, a few measurability conditions are required, as explained by Ernst et al. [3] in Subsection 2.3 of their work.

The existing PCE is largely founded on the independence assumption of input random variables. The assumption exploits product-type probability measures, enabling easy construction of multivariate orthogonal polynomials via tensorization of the spaces of univariate orthogonal polynomials. In reality, there may exist significant correlation or dependence among input variables, hindering or invalidating most existing stochastic methods, including PCE. For a general Gaussian input vector, there are at least two possibilities: (1) use a linear transformation to decorrelate the random variables and work with independent Gaussian variables; and (2) construct a sequence of weakly or strongly orthogonal multivariate polynomials consistent with the Gaussian measure and work with dependent Gaussian variables. Employing Gram-Schmidt orthogonalization [12], Navarro et al. [13] discussed construction of multivariate orthogonal polynomials for correlated variables. However, existence of Hermite orthogonal polynomials, which can be used as a basis for dependent Gaussian measures, has not been recognized. Soize and Ghanem [14] proposed orthogonal bases with regard to a general dependent probability measure of random input, but the bases are not necessarily polynomials. In consequence, analytical treatment of PCE-based statistics is highly non-trivial, if not impossible. In both works, a fundamental concern raised by Ernst et al. [3] about the completeness of orthogonal polynomials has not been addressed. Indeed, as demonstrated in this paper, the completeness is essential for the convergence of PCE subject to dependent Gaussian variables.

The main objective of this study is to generalize the classical PCE to account for arbitrary but dependent Gaussian probability measures. The paper is organized as follows. Section 2 defines or discusses mathematical notations and preliminaries. A brief exposition of multivariate orthogonal polynomials consistent with a general probability measure, including definitions of weak and strong orthogonalities, is given in Section 3. The section also describes relevant polynomial spaces and construction of their orthogonal decompositions. Section 4 defines multivariate Hermite polynomials consistent with a general dependent Gaussian probability measure. Two propositions proven herein reveal analytical formulae for the second-moment properties of these polynomials. The orthogonal basis and completeness of Hermite polynomials have also been discussed or proved. Section 5 formally presents a generalized PCE applicable for a general dependent Gaussian probability measure. The convergence, exactness, and optimality of the generalized PCE are explained. In the same section, the approximation quality of a truncated generalized PCE is discussed. The formulae for the mean and variance of the truncated generalized PCE are also derived. The application of the generalized PCE for infinitely many random variables is clarified. The section ends with a brief explanation on how and when the generalized PCE proposed can be extended for non-Gaussian probability measures. Numerical results from three illuminating examples, including a practical engineering problem, are reported in Section 6. Finally, conclusions are drawn in Section 7.

## 2. Notation and preliminaries

Let $\mathbb{N} := \{1, 2, \ldots\}$, $\mathbb{N}_0 := \mathbb{N} \cup \{0\}$, $\mathbb{R} := (-\infty, +\infty)$, $\mathbb{R}_0^+ := [0, +\infty)$, and $\mathbb{R}^+ := (0, +\infty)$ represent the sets of positive integer (natural), non-negative integer, real, non-negative real, and positive real numbers, respectively. For $N \in \mathbb{N}$, an $N$-dimensional multi-index is denoted by $\mathbf{j} := (j_1, \ldots, j_N) \in \mathbb{N}_0^N$ with degree $|\mathbf{j}| := j_1 + \cdots + j_N$ and factorial $\mathbf{j}! := j_1! \cdots j_N!$. These standard notations will be used throughout the paper.

Let $(\Omega, \mathcal{F}, \mathbb{P})$ be a complete probability space, where $\Omega$ is a sample space representing an abstract set of elementary events, $\mathcal{F}$ is a $\sigma$-field on $\Omega$, and $P : \mathcal{F} \to [0, 1]$ is a probability measure. With $\mathcal{B}^N$ representing the Borel $\sigma$-field on $\mathbb{R}^N$, $N \in \mathbb{N}$, consider an $\mathbb{R}^N$-valued Gaussian random vector $\mathbf{X} := (X_1, \ldots, X_N)^T : (\Omega, \mathcal{F}) \to (\mathbb{R}^N, \mathcal{B}^N)$, describing the statistical uncertainties in all system parameters of a stochastic problem. The input random variables are also referred to as basic random variables [3]. The non-zero, finite integer $N$ represents the number of input random variables and is referred to as the dimension of the stochastic problem.

Without loss of generality assume that $\mathbf{X}$ has a *zero* mean, that is, $\boldsymbol{\mu}_\mathbf{X} := \mathbb{E}[\mathbf{X}] = \mathbf{0} \in \mathbb{R}^N$; a symmetric, positive-definite covariance matrix $\boldsymbol{\Sigma}_\mathbf{X} := \mathbb{E}[\mathbf{X}\mathbf{X}^T] \in \mathbb{S}_+^N$, where $\mathbb{S}_+^N \subseteq \mathbb{R}^{N \times N}$ is the set of $N \times N$ real-valued, symmetric, positive-definite matrices; and a joint probability density function $\phi_\mathbf{X} : \mathbb{R}^N \to \mathbb{R}^+$, expressed by

$$\phi_\mathbf{X}(\mathbf{x}; \boldsymbol{\Sigma}_\mathbf{X}) := (2\pi)^{-\frac{N}{2}} (\det \boldsymbol{\Sigma}_\mathbf{X})^{-\frac{1}{2}} \exp\left[-\frac{1}{2}\mathbf{x}^T \boldsymbol{\Sigma}_\mathbf{X}^{-1} \mathbf{x}\right]. \tag{1}$$

Here, $\mathbb{E}$ is the expectation operator with respect to the probability measure $\mathbb{P}$ and $\det \boldsymbol{\Sigma}_\mathbf{X}$ is the determinant of $\boldsymbol{\Sigma}_\mathbf{X}$. Given the abstract probability space $(\Omega, \mathcal{F}, \mathbb{P})$ of $\mathbf{X}$, the image probability space is $(\mathbb{R}^N, \mathcal{B}^N, \phi_\mathbf{X} d\mathbf{x})$, where $\mathbb{R}^N$ can be viewed as the image of $\Omega$ from the mapping $\mathbf{X} : \Omega \to \mathbb{R}^N$, and is also the support of $\phi_{\mathbf{X}}(\mathbf{x}; \Sigma_{\mathbf{X}})$. The image probability space is convenient to use for computations. Indeed, relevant statements and objects in one space has obvious counterparts in the other space. Given any random variable $Y : (\Omega, \mathcal{F}) \to (\mathbb{R}, \mathbb{B})$, the Doob-Dynkin Lemma assures existence of a function $y : \mathbb{R}^N \to \mathbb{R}$ such that $Y(\omega) = y(\mathbf{X}(\omega))$. Furthermore, if $y$ is integrable, then the expectation of $Y$ can be defined by $\mathbb{E}[Y] := \int_\Omega Y(\omega) d\mathbb{P}(\omega)$ or $\mathbb{E}[Y] := \int_{\mathbb{R}^N} y(\mathbf{x}) \phi_{\mathbf{X}}(\mathbf{x}; \Sigma_{\mathbf{X}}) d\mathbf{x}$.

## 3. General multivariate orthogonal polynomials

For $\mathbf{j} \in \mathbb{N}_0^N$ and $\mathbf{x} = (x_1, \ldots, x_N) \in \mathbb{A}^N \subseteq \mathbb{R}^N$, a monomial in the variables $x_1, \ldots, x_N$ is the product $\mathbf{x}^\mathbf{j} = x_1^{j_1} \cdots x_N^{j_N}$ and has a total degree $|\mathbf{j}| = j_1 + \cdots + j_N$. A linear combination of $\mathbf{x}^\mathbf{j}$, where $|\mathbf{j}| = l \in \mathbb{N}_0$, is a homogeneous polynomial of degree $l$. Denote by
$$\mathcal{P}_l^N := \mathrm{span}\{\mathbf{x}^\mathbf{j} : |\mathbf{j}| = l, \mathbf{j} \in \mathbb{N}_0^N\}$$
the space of homogeneous polynomials of degree $l$, by
$$\Pi_m^N := \mathrm{span}\{\mathbf{x}^\mathbf{j} : 0 \le |\mathbf{j}| \le m, \mathbf{j} \in \mathbb{N}_0^N\}$$
the space of polynomials of degree at most $m \in \mathbb{N}_0$, and by
$$\Pi^N = \mathbb{R}[x_1, \ldots, x_N]$$
the space of all polynomials in $N$ variables. It is well known that the dimensions of the vector spaces $\mathcal{P}_l^N$ and $\Pi_m^N$, respectively, are [15]
$$\dim \mathcal{P}_l^N = \#\left\{\mathbf{j} \in \mathbb{N}_0^N : |\mathbf{j}| = l\right\} = \binom{N+l-1}{l}$$
and
$$\dim \Pi_m^N = \sum_{l=0}^m \dim \mathcal{P}_l^N = \sum_{l=0}^m \binom{N+l-1}{l} = \binom{N+m}{m}.$$

### 3.1. Measure-consistent orthogonal polynomials

Let $\mathbf{X} := (X_1, \ldots, X_N)^T$, $N \in \mathbb{N}$, be a general input random vector, which (1) has an absolutely continuous joint distribution function $F_{\mathbf{X}}(\mathbf{x})$ and a continuous joint probability density function $f_{\mathbf{X}}(\mathbf{x}) := \partial^N F_{\mathbf{X}}(\mathbf{x})/\partial x_1 \cdots \partial x_N$ with support $\mathbb{A}^N \subseteq \mathbb{R}^N$; and (2) possesses absolute finite moments of all orders, that is, for all $\mathbf{j} \in \mathbb{N}_0^N$,
$$\mu_\mathbf{j} := \mathbb{E}\left[|\mathbf{X}^\mathbf{j}|\right] := \int_{\mathbb{A}^N} |\mathbf{x}^\mathbf{j}| f_{\mathbf{X}}(\mathbf{x}) d\mathbf{x} < \infty.$$

For any polynomial pair $P, Q \in \Pi^N$, define an inner product
$$(P, Q)_{f_{\mathbf{X}} d\mathbf{x}} := \int_{\mathbb{A}^N} P(\mathbf{x}) Q(\mathbf{x}) f_{\mathbf{X}}(\mathbf{x}) d\mathbf{x} =: \mathbb{E}[P(\mathbf{X}) Q(\mathbf{X})] \tag{2}$$
with respect to the probability measure $f_{\mathbf{X}}(\mathbf{x}) d\mathbf{x}$ and the induced norm
$$\|P\|_{f_{\mathbf{X}} d\mathbf{x}} := \sqrt{(P, P)_{f_{\mathbf{X}} d\mathbf{x}}} = \left(\int_{\mathbb{A}^N} P^2(\mathbf{x}) f_{\mathbf{X}}(\mathbf{x}) d\mathbf{x}\right)^{1/2} = \sqrt{\mathbb{E}[P^2(\mathbf{X})]}.$$

The polynomials $P \in \Pi^N$ and $Q \in \Pi^N$ are called orthogonal to each other with respect to $f_{\mathbf{X}}(\mathbf{x}) d\mathbf{x}$ if $(P, Q)_{f_{\mathbf{X}} d\mathbf{x}} = 0$. This leads to a formal definition of multivariate orthogonal polynomials as follows.

**Definition 1** (Dunkl and Xu [15]). *A polynomial $P \in \Pi_l^N \subset \Pi^N$ is said to be an orthogonal polynomial of degree $l \in \mathbb{N}$ with respect to the inner product $(\cdot, \cdot)_{f_{\mathbf{X}} d\mathbf{x}}$, or alternatively with respect to the probability measure $f_{\mathbf{X}}(\mathbf{x}) d\mathbf{x}$, if it is orthogonal to all polynomials of lower degrees, that is, if*
$$(P, Q)_{f_{\mathbf{X}} d\mathbf{x}} = 0 \ \forall Q \in \Pi^N \text{ with } \deg Q < \deg P.$$

Under the prescribed assumptions, absolute moments of $\mathbf{X}$ of all orders exist, including the *zero*-order moment $\mu_0 := \int_{\mathbb{A}^N} f_{\mathbf{X}}(\mathbf{x})d\mathbf{x} = 1$ that is always positive. Evidently, $\|P\|_{f_{\mathbf{X}}d\mathbf{x}} > 0$ for all non-zero $P \in \Pi^N$. Then the inner product defined in (2) is positive-definite on $\Pi^N$. Therefore, there exists an infinite set of multivariate orthogonal polynomials [15], say, $\{P_{\mathbf{j}}(\mathbf{x}) : \mathbf{j} \in \mathbb{N}_0^N\}$, $P_{\mathbf{0}} = 1$, $P_{\mathbf{j}} \neq 0$, which is consistent with the probability measure $f_{\mathbf{X}}(\mathbf{x})d\mathbf{x}$, satisfying

$$\left(P_{\mathbf{j}}, P_{\mathbf{k}}\right)_{f_{\mathbf{X}}d\mathbf{x}} = 0 \text{ whenever } |\mathbf{j}| \neq |\mathbf{k}| \tag{3}$$

for $\mathbf{k} \in \mathbb{N}_0^N$. Here, the multi-index $\mathbf{j}$ of the multivariate polynomial $P_{\mathbf{j}}(\mathbf{x})$ refers to its total degree $|\mathbf{j}| = j_1 + \cdots + j_N$. Clearly, each $P_{\mathbf{j}} \in \Pi^N$ is an orthogonal polynomial according to Definition 1. This means that $P_{\mathbf{j}}$ is orthogonal to all polynomials of different degrees, but it may not be orthogonal to other orthogonal polynomials of the same degree.

Let $\mathcal{V}_0^N := \Pi_0^N = \text{span}\{1\}$ be the space of constant functions. For each $1 \leq l < \infty$, denote by $\mathcal{V}_l^N \subset \Pi_l^N$ the space of orthogonal polynomials of degree exactly $l$ that are orthogonal to all polynomials in $\Pi_{l-1}^N$, that is,

$$\mathcal{V}_l^N := \{P \in \Pi_l^N : (P, Q)_{f_{\mathbf{X}}d\mathbf{x}} = 0 \ \forall \ Q \in \Pi_{l-1}^N\}, \ 1 \leq l < \infty.$$

Then $\mathcal{V}_l^N$, provided that the support of $f_{\mathbf{X}}(\mathbf{x})$ has non-empty interior, is a vector space of dimension [15]

$$K_{N,l} := \dim \mathcal{V}_l^N = \dim \mathcal{P}_l^N = \binom{N + l - 1}{l}.$$

Many choices exist for the basis of $\mathcal{V}_l^N$; the bases of $\mathcal{V}_l^N$ do not have to be mutually orthogonal. With the exception of the monic orthogonal polynomials, the bases are not unique in the multivariate case. Here, to be formally proved in the next section, select $\{P_{\mathbf{j}}(\mathbf{x}) : |\mathbf{j}| = l, \mathbf{j} \in \mathbb{N}_0^N\} \subset \mathcal{V}_l^N$ to be a basis of $\mathcal{V}_l^N$, comprising $K_{N,l}$ number of basis functions. Each basis function $P_{\mathbf{j}}(\mathbf{x})$ is a multivariate orthogonal polynomial of degree $|\mathbf{j}|$ as discussed earlier. Obviously,

$$\mathcal{V}_l^N = \text{span}\{P_{\mathbf{j}} : |\mathbf{j}| = l, \mathbf{j} \in \mathbb{N}_0^N\}, \ 0 \leq l < \infty.$$

According to (3), $P_{\mathbf{j}}$ is orthogonal to $P_{\mathbf{k}}$ whenever $|\mathbf{j}| \neq |\mathbf{k}|$. Therefore, any two polynomial subspaces $\mathcal{V}_l^N$ and $\mathcal{V}_r^N$, where $0 \leq l, r < \infty$, are orthogonal whenever $l \neq r$. In consequence, there exist orthogonal decompositions of

$$\Pi_m^N = \bigoplus_{l=0}^m \mathcal{V}_l^N = \bigoplus_{l=0}^m \text{span}\{P_{\mathbf{j}} : |\mathbf{j}| = l, \mathbf{j} \in \mathbb{N}_0^N\} = \text{span}\{P_{\mathbf{j}} : 0 \leq |\mathbf{j}| \leq m, \mathbf{j} \in \mathbb{N}_0^N\}$$

and

$$\Pi^N = \bigoplus_{l \in \mathbb{N}_0} \mathcal{V}_l^N = \bigoplus_{l \in \mathbb{N}_0} \text{span}\{P_{\mathbf{j}} : |\mathbf{j}| = l, \mathbf{j} \in \mathbb{N}_0^N\} = \text{span}\{P_{\mathbf{j}} : \mathbf{j} \in \mathbb{N}_0^N\} \tag{4}$$

with the symbol $\oplus$ representing orthogonal sum.

### 3.2. Weak and strong orthogonalities

A possible lack of orthogonality between two distinct polynomials of the same degree can be used to characterize the strength of the orthogonality. Indeed, the multivariate orthogonal polynomials can be weakly orthogonal or strongly orthogonal.

**Definition 2.** *let $\mathbf{X} := (X_1, \ldots, X_N)^T$, $N \in \mathbb{N}$, be a general input random vector, which has a continuous joint probability density function $f_{\mathbf{X}}(\mathbf{x})$ with support $\mathbb{A}^N \subseteq \mathbb{R}^N$ and possesses absolute finite moments of all orders. Then a set of multivariate orthogonal polynomials $\{P_{\mathbf{j}}(\mathbf{x}) : \mathbf{j} \in \mathbb{N}_0^N\}$ consistent with the probability measure $f_{\mathbf{X}}(\mathbf{x})d\mathbf{x}$ is called a weakly orthogonal system if, for all $\mathbf{j}, \mathbf{k} \in \mathbb{N}_0^N$,*

$$\left(P_{\mathbf{j}}, P_{\mathbf{k}}\right)_{f_{\mathbf{X}}d\mathbf{x}} = 0 \text{ whenever } |\mathbf{j}| \neq |\mathbf{k}|;$$

and a strongly orthogonal system if, for all $\mathbf{j}, \mathbf{k} \in \mathbb{N}_0^N$,

$$\left(P_\mathbf{j}, P_\mathbf{k}\right)_{f_\mathbf{X} d\mathbf{x}} = 0 \text{ whenever } \mathbf{j} \neq \mathbf{k}.$$

Obviously, if a polynomial system is strongly orthogonal, then it is also weakly orthogonal. However, the converse is not true in general, for instance, when the variables are statistically dependent. When $\mathbf{X}$ comprises independent variables, then the multivariate polynomial system, if obtained via usual tensorized construction of the univariate polynomial spaces, becomes both strongly and weakly orthogonal. Nonetheless, Definition 2 can still be relevant for independent variables if the basis of $\mathcal{V}_l^N$ is chosen not to be orthogonal.

## 4. Multivariate Hermite orthogonal polynomials

When $\mathbf{X}$ has a Gaussian density function with support $\mathbb{R}^N$, as defined by (1), the moments $\int_{\mathbb{R}^N} |\mathbf{x}^\mathbf{j}| \phi_\mathbf{X}(\mathbf{x}; \Sigma_\mathbf{X}) d\mathbf{x}$ exist and are finite for all $\mathbf{j} \in \mathbb{N}_0^N$. Therefore, orthogonal polynomials in $\mathbf{x}$ exist with respect to the inner product

$$(P, Q)_{\phi_\mathbf{X} d\mathbf{x}} := \int_{\mathbb{R}^N} P(\mathbf{x}) Q(\mathbf{x}) \phi_\mathbf{X}(\mathbf{x}; \Sigma_\mathbf{X}) d\mathbf{x} =: \mathbb{E}\left[P(\mathbf{X}) Q(\mathbf{X})\right] \tag{5}$$

or the probability measure $\phi_\mathbf{X}(\mathbf{x}; \Sigma_\mathbf{X}) d\mathbf{x}$. Here, a special basis of $\mathcal{V}_l^N$, denoted by $\{H_\mathbf{j}(\mathbf{x}; \Sigma_\mathbf{X}) : |\mathbf{j}| = l, \mathbf{j} \in \mathbb{N}_0^N\} \subset \mathcal{V}_l^N$, is presented, which will be proved later to be weakly orthogonal as per Definition 2. The set of all such polynomials, that is, $\{H_\mathbf{j}(\mathbf{x}; \Sigma_\mathbf{X}) : \mathbf{j} \in \mathbb{N}_0^N\} \subset \Pi^N$ comprises polynomials that are orthogonal with respect to the inner product in (5). The polynomials are consistent with the probability measure $\phi_\mathbf{X}(\mathbf{x}; \Sigma_\mathbf{X}) d\mathbf{x}$ and are often referred to as multivariate Hermite orthogonal polynomials.

### 4.1. Definition

A popular approach for defining multivariate Hermite polynomials entails derivatives of the multivariate Gaussian probability density function. Many researchers have used this definition [16, 17, 18, 19]. Formal definitions of both orthogonal and standardized orthogonal polynomials follow.

**Definition 3.** *Let $\mathbf{X} = (X_1, \ldots, X_N)^T$, $N \in \mathbb{N}$, be an $\mathbb{R}^N$-valued Gaussian random vector with zero mean; symmetric, positive-definite covariance matrix $\Sigma_\mathbf{X} \in S_+^N$; and multivariate density function $\phi_\mathbf{X}(\mathbf{x}; \Sigma_\mathbf{X})$. Then a multivariate Hermite orthogonal polynomial $H_\mathbf{j}(\mathbf{x}; \Sigma_\mathbf{X})$, $\mathbf{j} = (j_1, \ldots, j_N) \in \mathbb{N}_0^N$, of degree $|\mathbf{j}| = j_1 + \cdots + j_N$ is defined as*

$$H_\mathbf{j}(\mathbf{x}; \Sigma_\mathbf{X}) := \frac{(-1)^{|\mathbf{j}|}}{\phi_\mathbf{X}(\mathbf{x}; \Sigma_\mathbf{X})} \left(\frac{\partial}{\partial \mathbf{x}}\right)^\mathbf{j} \phi_\mathbf{X}(\mathbf{x}; \Sigma_\mathbf{X}), \tag{6}$$

*where $(\partial/\partial \mathbf{x})^\mathbf{j} := \partial^{j_1+\cdots+j_N}/\partial x_1^{j_1} \cdots \partial x_N^{j_N}$.*

**Definition 4.** *A standardized multivariate Hermite orthogonal polynomial $\Psi_\mathbf{j}(\mathbf{x}; \Sigma_\mathbf{X})$, $\mathbf{j} = (j_1, \ldots, j_N) \in \mathbb{N}_0^N$, of degree $|\mathbf{j}| = j_1 + \cdots + j_N$ is defined as*

$$\Psi_\mathbf{j}(\mathbf{x}; \Sigma_\mathbf{X}) := \frac{H_\mathbf{j}(\mathbf{x}; \Sigma_\mathbf{X})}{(H_\mathbf{j}(\mathbf{x}; \Sigma_\mathbf{X}), H_\mathbf{j}(\mathbf{x}; \Sigma_\mathbf{X}))_{\phi_\mathbf{X} d\mathbf{x}}} = \frac{H_\mathbf{j}(\mathbf{x}; \Sigma_\mathbf{X})}{\sqrt{\mathbb{E}[H_\mathbf{j}^2(\mathbf{X}; \Sigma_\mathbf{X})]}}. \tag{7}$$

Definition 3 is a generalization of the definition of the $j$th-degree univariate Hermite orthogonal polynomial

$$H_j(x) = \frac{(-1)^j}{\phi_X(x)} \frac{d^j}{dx^j} \phi_X(x), \quad j \in \mathbb{N}_0,$$

known as Rodrigues's formula [20], where $\phi_X(x) = (2\pi)^{-1/2} \exp(-x^2/2)$ is the probability density function of a *zero*-mean Gaussian random variable with unit variance. Definition 4 facilitates scaling of multivariate Hermite polynomials, so that their standardized version reduces to multivariate orthonormal polynomials for independent random variables. The standardized multivariate polynomials should not be confused with multivariate orthonormal polynomials for dependent random variables.

If the Gaussian random variables are independent, then the covariance matrix becomes diagonal, that is, $\Sigma_\mathbf{X} = \text{diag}(\sigma_1^2, \ldots, \sigma_N^2)$ with $0 < \sigma_i^2 < \infty$, $i = 1, \ldots, N$, representing the variance of the $i$th variable. If, in addition, $\sigma_i^2 = 1$ for all $i = 1, \ldots, N$, then $\Sigma_\mathbf{X} = \mathbf{I}$, the $N$-dimensional identity matrix, and (1) leads to a product-type density function $\phi_\mathbf{X}(\mathbf{x}; \mathbf{I}) = \Pi_{i=1}^N \phi_{X_i}(x_i)$, comprising marginal probability density functions $\phi_{X_i}(x_i) = (2\pi)^{-1/2} \exp(-x_i^2/2)$, $i = 1, \ldots, N$. In consequence, Definition 3 simplifies to the well-known tensorized construction: $H_\mathbf{j}(\mathbf{x}; \mathbf{I}) = H_{j_1}(x_1) \cdots H_{j_N}(x_N)$, that is, a multivariate orthogonal polynomial of degree $|\mathbf{j}|$ is simply a product of $N$ univariate orthogonal polynomials $H_{j_i}(x_i)$, $i = 1, \ldots, N$, of degree $j_i$ such that $j_1 + \cdots + j_N = |\mathbf{j}|$.

According to Definition 3, the set of Hermite polynomials $\{H_\mathbf{j}(\mathbf{x}; \Sigma_\mathbf{X}), \mathbf{j} \in \mathbb{N}_0^N\}$ for general dependent Gaussian variables is weakly orthogonal with respect to $(\cdot, \cdot)_{\phi_\mathbf{X} d\mathbf{x}}$, that is,

$$\left(H_\mathbf{j}, H_\mathbf{k}\right)_{\phi_\mathbf{X} d\mathbf{x}} = \mathbb{E}\left[H_\mathbf{j}(\mathbf{X}; \Sigma_\mathbf{X}) H_\mathbf{k}(\mathbf{X}; \Sigma_\mathbf{X})\right] = 0, \ |\mathbf{j}| \neq |\mathbf{k}|,$$

to be formally proved in the following subsection. This means that $H_\mathbf{j}$ is orthogonal to all polynomials of different degrees, but it may not be orthogonal to other orthogonal polynomials of the same degree. However, if the Gaussian variables are independent, then the resultant multivariate Hermite polynomials are strongly orthogonal. This is because of the product structure of such polynomials, where any two univariate Hermite polynomials of distinct degrees are orthogonal. Since the focus of this work is dependent Gaussian variables, the orthogonality of multivariate Hermite polynomials for the rest of the paper should be interpreted in the context of weak orthogonality.

### 4.2. Second-moment properties

When the input random variables $X_1, \ldots, X_N$, instead of the variables $x_1, \ldots, x_N$, are inserted in the argument, the Hermite orthogonal polynomials become random functions of Gaussian input vector $\mathbf{X} = (X_1, \ldots, X_N)^T$. Therefore, it is important to derive explicit formulae for their second-moment properties in terms of the statistics of $\mathbf{X}$. The formulae, obtained here using a compact form of the generating function in Proposition 6, are described by Propositions 7 and 8.

**Definition 5.** *The generating function for the family of multivariate Hermite orthogonal polynomials $\{H_\mathbf{j}(\mathbf{x}; \Sigma_\mathbf{X}), \mathbf{j} \in \mathbb{N}_0^N\}$ is defined as the convergent expansion*

$$\sum_{\mathbf{j} \in \mathbb{N}_0^N} \frac{\mathbf{t}^\mathbf{j}}{\mathbf{j}!} H_\mathbf{j}(\mathbf{x}; \Sigma_\mathbf{X}), \ \mathbf{t} \in \mathbb{R}^N, \tag{8}$$

*where $\mathbf{j}! := j_1! \ldots j_N!$ and $\mathbf{t}^\mathbf{j} := t_1^{j_1} \ldots t_N^{j_N}$.*

**Proposition 6.** *In reference to Definition 5, the generating function for $\mathbf{t} \in \mathbb{R}^N$ is*

$$\sum_{\mathbf{j} \in \mathbb{N}_0^N} \frac{\mathbf{t}^\mathbf{j}}{\mathbf{j}!} H_\mathbf{j}(\mathbf{x}; \Sigma_\mathbf{X}) = \exp\left(\mathbf{t}^T \Sigma_\mathbf{X}^{-1} \mathbf{x} - \frac{1}{2} \mathbf{t}^T \Sigma_\mathbf{X}^{-1} \mathbf{t}\right), \tag{9}$$

*where $\Sigma_\mathbf{X}^{-1}$ is the inverse of $\Sigma_\mathbf{X}$ and the symbol $T$ denotes matrix transposition.*

*Proof.* Using the definition of $H_\mathbf{j}(\mathbf{x}; \Sigma_\mathbf{X})$ from (6),

$$\begin{aligned}
\sum_{\mathbf{j} \in \mathbb{N}_0^N} \frac{\mathbf{t}^\mathbf{j}}{\mathbf{j}!} H_\mathbf{j}(\mathbf{x}; \Sigma_\mathbf{X}) &= \frac{1}{\phi_\mathbf{X}(\mathbf{x}; \Sigma_\mathbf{X})} \sum_{\mathbf{j} \in \mathbb{N}_0^N} \frac{\mathbf{t}^\mathbf{j}}{\mathbf{j}!} (-1)^{|\mathbf{j}|} \left(\frac{\partial}{\partial \mathbf{x}}\right)^\mathbf{j} \phi_\mathbf{X}(\mathbf{x}; \Sigma_\mathbf{X}) \\
&= \frac{\phi_\mathbf{X}(\mathbf{x} - \mathbf{t}; \Sigma_\mathbf{X})}{\phi_\mathbf{X}(\mathbf{x}; \Sigma_\mathbf{X})} \\
&= \exp\left(\mathbf{t}^T \Sigma_\mathbf{X}^{-1} \mathbf{x} - \frac{1}{2} \mathbf{t}^T \Sigma_\mathbf{X}^{-1} \mathbf{t}\right).
\end{aligned}$$

Here, the second line is formed by recognizing the sum in the first equality to be the Taylor series expansion of $\phi_{\mathbf{X}}(\mathbf{x} - \mathbf{t}; \Sigma_{\mathbf{X}})$ at $\mathbf{x}$, whereas the third line is obtained by applying (1) and reduction. $\square$

**Proposition 7.** *The first-order moments of multivariate Hermite orthogonal polynomials are*

$$\mathbb{E}\left[H_{\mathbf{j}}(\mathbf{X}; \Sigma_{\mathbf{X}})\right] = \begin{cases} 1, & \mathbf{j} = \mathbf{0}, \\ 0, & \mathbf{j} \neq \mathbf{0}. \end{cases} \tag{10}$$

*Proof.* Multiplying the generating function for $\mathbf{t} \in \mathbb{R}^N$ in (8) with $\phi_{\mathbf{X}}(\mathbf{x}; \Sigma_{\mathbf{X}})$ and then integrating over $\mathbb{R}^N$ gives

$$\begin{aligned}
& \int_{\mathbb{R}^N} \sum_{\mathbf{j} \in \mathbb{N}_0^N} \frac{\mathbf{t}^{\mathbf{j}}}{\mathbf{j}!} H_{\mathbf{j}}(\mathbf{x}; \Sigma_{\mathbf{X}}) \phi_{\mathbf{X}}(\mathbf{x}; \Sigma_{\mathbf{X}}) d\mathbf{x} \\
& = \int_{\mathbb{R}^N} \exp\left(\mathbf{t}^T \Sigma_{\mathbf{X}}^{-1} \mathbf{x} - \frac{1}{2} \mathbf{t}^T \Sigma_{\mathbf{X}}^{-1} \mathbf{t}\right) (2\pi)^{-\frac{N}{2}} (\det \Sigma_{\mathbf{X}})^{-\frac{1}{2}} \exp\left(-\frac{1}{2} \mathbf{x}^T \Sigma_{\mathbf{X}}^{-1} \mathbf{x}\right) d\mathbf{x} \\
& = \int_{\mathbb{R}^N} (2\pi)^{-\frac{N}{2}} (\det \Sigma_{\mathbf{X}})^{-\frac{1}{2}} \exp\left\{-\frac{1}{2}(\mathbf{x} - \mathbf{t})^T \Sigma_{\mathbf{X}}^{-1}(\mathbf{x} - \mathbf{t})\right\} d\mathbf{x} \\
& = 1,
\end{aligned} \tag{11}$$

where the second line uses Proposition 6, that is, (9), and (1); the third line is obtained by reduction, yielding unity in the last line – the result of integrating a Gaussian probability density function on $\mathbb{R}^N$. Finally, comparing the coefficients of $\mathbf{t}^{\mathbf{j}}$, $\mathbf{j} \in \mathbb{N}_0^N$, in (11) produces

$$\int_{\mathbb{R}^N} H_{\mathbf{j}}(\mathbf{x}; \Sigma_{\mathbf{X}}) \phi_{\mathbf{X}}(\mathbf{x}; \Sigma_{\mathbf{X}}) d\mathbf{x} = \begin{cases} 1, & \mathbf{j} = \mathbf{0}, \\ 0, & \mathbf{j} \neq \mathbf{0}, \end{cases}$$

where the integral on the left is the same as the first-order moment, hence completing the proof. $\square$

**Proposition 8.** *The second-order moments of multivariate Hermite orthogonal polynomials are*

$$\mathbb{E}\left[H_{\mathbf{j}}(\mathbf{X}; \Sigma_{\mathbf{X}}) H_{\mathbf{k}}(\mathbf{X}; \Sigma_{\mathbf{X}})\right] = \begin{cases} \mathbf{j}! \mathbf{k}! \sum_{\substack{\boldsymbol{\theta} \in \mathbb{N}_0^{N \times N} \\ \mathbf{r}(\boldsymbol{\theta}) = \mathbf{j}, \mathbf{c}(\boldsymbol{\theta}) = \mathbf{k} \\ |\mathbf{j}| = |\mathbf{k}|}} \frac{\left(\Sigma_{\mathbf{X}}^{-1}\right)^{\boldsymbol{\theta}}}{\boldsymbol{\theta}!}, & |\mathbf{j}| = |\mathbf{k}|, \\ 0, & |\mathbf{j}| \neq |\mathbf{k}|, \end{cases} \tag{12}$$

*where $\boldsymbol{\theta} \in \mathbb{N}_0^{N \times N}$ is an index matrix, comprising non-negative integers, with the $(p,q)$th element $\theta_{pq} \in \mathbb{N}_0$ for $p, q = 1, \ldots, N$; $\mathbf{r}(\boldsymbol{\theta}) = (r_1, \ldots, r_N)$ is the row-sum vector of $\boldsymbol{\theta}$ with the $p$th element $r_p = \sum_{q=1}^N \theta_{pq}$; $\mathbf{c}(\boldsymbol{\theta}) = (c_1, \ldots, c_N)$ is the column-sum vector of $\boldsymbol{\theta}$ with the $q$th element $c_q = \sum_{p=1}^N \theta_{pq}$;*

$$\boldsymbol{\theta}! := \prod_{p,q=1}^N \theta_{pq}!;$$

*and*

$$\left(\Sigma_{\mathbf{X}}^{-1}\right)^{\boldsymbol{\theta}} := \prod_{p,q=1}^N \left(\Sigma_{\mathbf{X},pq}^{-1}\right)^{\theta_{pq}}$$

*with $\Sigma_{\mathbf{X},pq}^{-1}$ representing the $(p,q)$th element of $\Sigma_{\mathbf{X}}^{-1}$. The summation in (12) is over all index matrices $\boldsymbol{\theta}$ with the* *row-sum vector* $\mathbf{r}(\boldsymbol{\theta}) = \mathbf{j}$ *and the column-sum vector* $\mathbf{c}(\boldsymbol{\theta}) = \mathbf{k}$ *such that* $|\mathbf{j}| = |\mathbf{k}|$. *Furthermore,*

$$\mathbb{E}\left[H_{\mathbf{j}}^2(\mathbf{X};\boldsymbol{\Sigma}_{\mathbf{X}})\right] = (\mathbf{j}!)^2 \sum_{\substack{\boldsymbol{\theta} \in \mathbb{N}_0^{N \times N} \\ \mathbf{r}(\boldsymbol{\theta})=\mathbf{c}(\boldsymbol{\theta})=\mathbf{j}}} \frac{\left(\boldsymbol{\Sigma}_{\mathbf{X}}^{-1}\right)^{\boldsymbol{\theta}}}{\boldsymbol{\theta}!}. \tag{13}$$

*Proof.* Multiplying the product of two generating functions for $\mathbf{t}, \mathbf{s} \in \mathbb{R}^N$ in (8) with $\phi_{\mathbf{X}}(\mathbf{x};\boldsymbol{\Sigma}_{\mathbf{X}})$ and then integrating over $\mathbb{R}^N$ gives

$$\begin{aligned}
&\int_{\mathbb{R}^N} \sum_{\mathbf{j} \in \mathbb{N}_0^N} \sum_{\mathbf{k} \in \mathbb{N}_0^N} \frac{\mathbf{t}^{\mathbf{j}}}{\mathbf{j}!} H_{\mathbf{j}}(\mathbf{x};\boldsymbol{\Sigma}_{\mathbf{X}}) \frac{\mathbf{s}^{\mathbf{k}}}{\mathbf{k}!} H_{\mathbf{k}}(\mathbf{x};\boldsymbol{\Sigma}_{\mathbf{X}}) \phi_{\mathbf{X}}(\mathbf{x};\boldsymbol{\Sigma}_{\mathbf{X}}) d\mathbf{x} \\
&= \int_{\mathbb{R}^N} \exp\left(\mathbf{t}^T \boldsymbol{\Sigma}_{\mathbf{X}}^{-1} \mathbf{x} - \frac{1}{2}\mathbf{t}^T \boldsymbol{\Sigma}_{\mathbf{X}}^{-1} \mathbf{t}\right) \exp\left(\mathbf{s}^T \boldsymbol{\Sigma}_{\mathbf{X}}^{-1} \mathbf{x} - \frac{1}{2}\mathbf{s}^T \boldsymbol{\Sigma}_{\mathbf{X}}^{-1} \mathbf{s}\right) \frac{\exp\left(-\frac{1}{2}\mathbf{x}^T \boldsymbol{\Sigma}_{\mathbf{X}}^{-1} \mathbf{x}\right)}{(2\pi)^{\frac{N}{2}} (\det \boldsymbol{\Sigma}_{\mathbf{X}})^{\frac{1}{2}}} d\mathbf{x} \\
&= \exp\left(\mathbf{t}^T \boldsymbol{\Sigma}_{\mathbf{X}}^{-1} \mathbf{s}\right) \int_{\mathbb{R}^N} (2\pi)^{-\frac{N}{2}} (\det \boldsymbol{\Sigma}_{\mathbf{X}})^{-\frac{1}{2}} \exp\left\{-\frac{1}{2}(\mathbf{x}-\mathbf{t}-\mathbf{s})^T \boldsymbol{\Sigma}_{\mathbf{X}}^{-1}(\mathbf{x}-\mathbf{t}-\mathbf{s})\right\} d\mathbf{x} \\
&= \exp\left(\mathbf{t}^T \boldsymbol{\Sigma}_{\mathbf{X}}^{-1} \mathbf{s}\right),
\end{aligned} \tag{14}$$

where the last line is obtained by applying again Proposition 6, that is, (9), using (1), and finally recognizing the Gaussian integral to be unity.

For $\boldsymbol{\Sigma}_{\mathbf{X}}^{-1} \in \mathbb{R}^{N \times N}$, one has the convergent expansion [21]

$$\exp\left(\mathbf{t}^T \boldsymbol{\Sigma}_{\mathbf{X}}^{-1} \mathbf{s}\right) = \sum_{\boldsymbol{\theta} \in \mathbb{N}_0^{N \times N}} \frac{\left(\boldsymbol{\Sigma}_{\mathbf{X}}^{-1}\right)^{\boldsymbol{\theta}}}{\boldsymbol{\theta}!} \mathbf{t}^{\mathbf{r}(\boldsymbol{\theta})} \mathbf{s}^{\mathbf{c}(\boldsymbol{\theta})}. \tag{15}$$

Let $\mathbf{j} = \mathbf{r}(\boldsymbol{\theta})$ and $\mathbf{k} = \mathbf{c}(\boldsymbol{\theta})$ in (15) and note that $|\mathbf{j}| = |\mathbf{k}|$. Therefore,

$$\exp\left(\mathbf{t}^T \boldsymbol{\Sigma}_{\mathbf{X}}^{-1} \mathbf{s}\right) = \sum_{l=0}^{\infty} \sum_{\substack{\mathbf{j} \in \mathbb{N}_0^N, \mathbf{k} \in \mathbb{N}_0^N \\ |\mathbf{j}|=|\mathbf{k}|=l}} \sum_{\substack{\boldsymbol{\theta} \in \mathbb{N}_0^{N \times N} \\ \mathbf{r}(\boldsymbol{\theta})=\mathbf{j}, \mathbf{c}(\boldsymbol{\theta})=\mathbf{k} \\ |\mathbf{j}|=|\mathbf{k}|=l}} \frac{\left(\boldsymbol{\Sigma}_{\mathbf{X}}^{-1}\right)^{\boldsymbol{\theta}}}{\boldsymbol{\theta}!} \mathbf{t}^{\mathbf{j}} \mathbf{s}^{\mathbf{k}}. \tag{16}$$

Interchanging the integral and summation operators of (14) and using (16) gives

$$\sum_{\mathbf{j} \in \mathbb{N}_0^N, \mathbf{k} \in \mathbb{N}_0^N} \int_{\mathbb{R}^N} \frac{\mathbf{t}^{\mathbf{j}}}{\mathbf{j}!} H_{\mathbf{j}}(\mathbf{x};\boldsymbol{\Sigma}_{\mathbf{X}}) \frac{\mathbf{s}^{\mathbf{k}}}{\mathbf{k}!} H_{\mathbf{k}}(\mathbf{x};\boldsymbol{\Sigma}_{\mathbf{X}}) \phi_{\mathbf{X}}(\mathbf{x};\boldsymbol{\Sigma}_{\mathbf{X}}) d\mathbf{x} = \sum_{l=0}^{\infty} \sum_{\substack{\mathbf{j} \in \mathbb{N}_0^N, \mathbf{k} \in \mathbb{N}_0^N \\ |\mathbf{j}|=|\mathbf{k}|=l}} \sum_{\substack{\boldsymbol{\theta} \in \mathbb{N}_0^{N \times N} \\ \mathbf{r}(\boldsymbol{\theta})=\mathbf{j}, \mathbf{c}(\boldsymbol{\theta})=\mathbf{k} \\ |\mathbf{j}|=|\mathbf{k}|=l}} \frac{\left(\boldsymbol{\Sigma}_{\mathbf{X}}^{-1}\right)^{\boldsymbol{\theta}}}{\boldsymbol{\theta}!} \mathbf{t}^{\mathbf{j}} \mathbf{s}^{\mathbf{k}}. \tag{17}$$

Finally, comparing the coefficients of $\mathbf{t}^{\mathbf{j}} \mathbf{s}^{\mathbf{k}}$, $\mathbf{j}, \mathbf{k} \in \mathbb{N}_0^N$, in (17) yields

$$\int_{\mathbb{R}^N} H_{\mathbf{j}}(\mathbf{x};\boldsymbol{\Sigma}_{\mathbf{X}}) H_{\mathbf{k}}(\mathbf{x};\boldsymbol{\Sigma}_{\mathbf{X}}) \phi_{\mathbf{X}}(\mathbf{x};\boldsymbol{\Sigma}_{\mathbf{X}}) d\mathbf{x} = \begin{cases} \mathbf{j}! \mathbf{k}! \displaystyle\sum_{\substack{\boldsymbol{\theta} \in \mathbb{N}_0^{N \times N} \\ \mathbf{r}(\boldsymbol{\theta})=\mathbf{j}, \mathbf{c}(\boldsymbol{\theta})=\mathbf{k} \\ |\mathbf{j}|=|\mathbf{k}|}} \dfrac{\left(\boldsymbol{\Sigma}_{\mathbf{X}}^{-1}\right)^{\boldsymbol{\theta}}}{\boldsymbol{\theta}!}, & |\mathbf{j}| = |\mathbf{k}|, \\ 0, & |\mathbf{j}| \neq |\mathbf{k}|, \end{cases} \tag{18}$$

where the integral on the left is the same as the second-order moment, hence obtaining the desired result in (12). Setting $\mathbf{j} = \mathbf{k}$ in (18) produces (13). □

From (6), the *zero*-degree Hermite orthogonal polynomial $H_{\mathbf{0}}(\mathbf{x};\boldsymbol{\Sigma}_{\mathbf{X}}) = 1$, regardless of $\boldsymbol{\Sigma}_{\mathbf{X}}$. For any $\mathbf{j} \in \mathbb{N}_0^N$ and $\mathbf{k} = \mathbf{0}$, (12) reproduces (10), the first-order moment of $H_\mathbf{j}(\mathbf{X}; \Sigma_\mathbf{X})$. Therefore, Proposition 8 subsumes Proposition 7.

**Corollary 9.** *The first- and second-order moments of standardized multivariate Hermite orthogonal polynomials are*

$$\mathbb{E}\left[\Psi_\mathbf{j}(\mathbf{X}; \Sigma_\mathbf{X})\right] = \begin{cases} 1, & \mathbf{j} = \mathbf{0}, \\ 0, & \mathbf{j} \neq \mathbf{0}, \end{cases}$$

*and*

$$\mathbb{E}\left[\Psi_\mathbf{j}(\mathbf{X}; \Sigma_\mathbf{X})\Psi_\mathbf{k}(\mathbf{X}; \Sigma_\mathbf{X})\right] = \begin{cases} \dfrac{\mathbb{E}\left[H_\mathbf{j}(\mathbf{X}; \Sigma_\mathbf{X})H_\mathbf{k}(\mathbf{X}; \Sigma_\mathbf{X})\right]}{\sqrt{\mathbb{E}[H_\mathbf{j}^2(\mathbf{X}; \Sigma_\mathbf{X})]}\sqrt{\mathbb{E}[H_\mathbf{k}^2(\mathbf{X}; \Sigma_\mathbf{X})]}}, & |\mathbf{j}| = |\mathbf{k}|, \\ 0, & |\mathbf{j}| \neq |\mathbf{k}|, \end{cases} \quad (19)$$

*respectively, including*

$$\mathbb{E}\left[\Psi_\mathbf{j}^2(\mathbf{X}; \Sigma_\mathbf{X})\right] = 1, \ \mathbf{j} \in \mathbb{N}_0^N, \quad (20)$$

*where the expectations*, $\mathbb{E}[H_\mathbf{j}(\mathbf{X}; \Sigma_\mathbf{X})H_\mathbf{k}(\mathbf{X}; \Sigma_\mathbf{X})]$ *and* $\mathbb{E}[H_\mathbf{j}^2(\mathbf{X}; \Sigma_\mathbf{X})]$, *are obtained from* (12) *and* (13), *respectively.*

**Corollary 10.** *If* $\mathbf{X} = (X_1, \ldots, X_N)^T$ *comprises independent Gaussian random variables, each with zero mean and unit variance, then* $\Sigma_\mathbf{X} = \mathbf{I}$, *the N-dimensional identity matrix, resulting in multivariate Hermite orthonormal polynomials* $[\Psi_\mathbf{j}(\mathbf{X}; \mathbf{I})]$, $\mathbf{j} \in \mathbb{N}_0^N$, *with their first- and second-order moments*

$$\mathbb{E}\left[\Psi_\mathbf{j}(\mathbf{X}; \mathbf{I})\right] = \begin{cases} 1, & \mathbf{j} = \mathbf{0}, \\ 0, & \mathbf{j} \neq \mathbf{0}, \end{cases}$$

*and*

$$\mathbb{E}\left[\Psi_\mathbf{j}(\mathbf{X}; \mathbf{I})\Psi_\mathbf{k}(\mathbf{X}; \mathbf{I})\right] = \begin{cases} 1, & \mathbf{j} = \mathbf{k}, \\ 0, & \mathbf{j} \neq \mathbf{k}, \end{cases}$$

*respectively, including*

$$\mathbb{E}\left[\Psi_\mathbf{j}^2(\mathbf{X}; \mathbf{I})\right] = 1, \ \mathbf{j} \in \mathbb{N}_0^N.$$

From Corollaries 9 and 10, the first-order moments of standardized Hermite orthogonal polynomials for dependent variables and Hermite orthonormal polynomials for independent variables are the same. However, the second-order moments of Hermite orthonormal polynomials for independent variables simplify significantly due to strong orthogonality. This explains why the development of the classical PCE for independent variables is not unduly difficult.

### 4.3. Orthogonal basis and completeness

An important question regarding Hermite orthogonal polynomials is whether they constitute a basis in a function space of interest, such as a Hilbert space. Let $L^2(\mathbb{R}^N, \mathcal{B}^N, \phi_\mathbf{X} d\mathbf{x})$ represent a Hilbert space of square-integrable functions with respect to the Gaussian probability measure $\phi_\mathbf{X}(\mathbf{x}; \Sigma_\mathbf{X}) d\mathbf{x}$ supported on $\mathbb{R}^N$. The following two propositions show that, indeed, Hermite orthogonal polynomials span various spaces of interest.

**Proposition 11.** *Let* $\mathbf{X} := (X_1, \ldots, X_N)^T : (\Omega, \mathcal{F}) \to (\mathbb{R}^N, \mathcal{B}^N)$, $N \in \mathbb{N}$, *be an* $\mathbb{R}^N$-*valued Gaussian random vector with zero mean; symmetric, positive-definite covariance matrix* $\Sigma_\mathbf{X}$; *and multivariate probability density function* $\phi_\mathbf{X}(\mathbf{x}; \Sigma_\mathbf{X})$. *Then,* $\{H_\mathbf{j}(\mathbf{x}; \Sigma_\mathbf{X}) : |\mathbf{j}| = l, \mathbf{j} \in \mathbb{N}_0^N\}$, *the set of multivariate Hermite orthogonal polynomials of degree l consistent with the Gaussian probability measure* $\phi_\mathbf{X} d\mathbf{x}$, *is a basis of* $\mathcal{V}_l^N$.

*Proof.* According to Takemura and Takeuchi [18], the multivariate Hermite polynomials from Definition 3 are orthogonal to their dual polynomials [18, 19]

$$\tilde{H}_\mathbf{j}(\mathbf{x}; \Sigma_\mathbf{X}) := \frac{(-1)^{|\mathbf{j}|}}{\phi_\mathbf{X}(\mathbf{x}; \Sigma_\mathbf{X})} \left(\frac{\partial}{\partial \mathbf{z}}\right)^\mathbf{j} \phi_\mathbf{X}(\Sigma_\mathbf{X} \mathbf{z}; \Sigma_\mathbf{X}), \ \mathbf{z} = \Sigma_\mathbf{X}^{-1} \mathbf{x},$$

in the sense that

$$\mathbb{E}\left[\tilde{H}_\mathbf{j}(\mathbf{X};\boldsymbol{\Sigma}_\mathbf{X})H_\mathbf{k}(\mathbf{X};\boldsymbol{\Sigma}_\mathbf{X})\right] := \int_{\mathbb{R}^N} \tilde{H}_\mathbf{j}(\mathbf{x};\boldsymbol{\Sigma}_\mathbf{X})H_\mathbf{k}(\mathbf{x};\boldsymbol{\Sigma}_\mathbf{X})\phi_\mathbf{X}(\mathbf{x};\boldsymbol{\Sigma}_\mathbf{X})d\mathbf{x} = \begin{cases} \mathbf{j}!, & \mathbf{j} = \mathbf{k}, \\ 0, & \mathbf{j} \neq \mathbf{k}. \end{cases} \quad (21)$$

Denote by $\mathbf{H}_l(\mathbf{x};\boldsymbol{\Sigma}_\mathbf{X}) = (H_{l,1}(\mathbf{x};\boldsymbol{\Sigma}_\mathbf{X}),\ldots,H_{l,K_{N,l}}(\mathbf{x};\boldsymbol{\Sigma}_\mathbf{X}))^T$ and $\tilde{\mathbf{H}}_l(\mathbf{x};\boldsymbol{\Sigma}_\mathbf{X}) = (\tilde{H}_{l,1}(\mathbf{x};\boldsymbol{\Sigma}_\mathbf{X}),\ldots,\tilde{H}_{l,K_{N,l}}(\mathbf{x};\boldsymbol{\Sigma}_\mathbf{X}))^T$ the column vectors of the elements of $\{H_\mathbf{j}(\mathbf{x};\boldsymbol{\Sigma}_\mathbf{X}) : |\mathbf{j}| = l, \mathbf{j} \in \mathbb{N}_0^N\}$ and $\{\tilde{H}_\mathbf{j}(\mathbf{x};\boldsymbol{\Sigma}_\mathbf{X}) : |\mathbf{j}| = l, \mathbf{j} \in \mathbb{N}_0^N\}$, respectively, both arranged according to some monomial order of choice. Let $\mathbf{a}_l^T = (a_{l,1},\ldots,a_{l,K_{N,l}})$ be a row vector comprising some constants $a_{l,i} \in \mathbb{R}$, $i = 1,\ldots,K_{N,l}$. Set $\mathbf{a}_l^T \mathbf{H}_l(\mathbf{x};\boldsymbol{\Sigma}_\mathbf{X}) = 0$. Multiply both sides of the equality from the right by $\tilde{\mathbf{H}}_l^T(\mathbf{x};\boldsymbol{\Sigma}_\mathbf{X})$, integrate with respect to the measure $\phi_\mathbf{X}(\mathbf{x};\boldsymbol{\Sigma}_\mathbf{X})d\mathbf{x}$ over $\mathbb{R}^N$, and apply transposition to obtain

$$\mathbf{G}_l \mathbf{a}_l = \mathbf{0}, \quad (22)$$

where $\mathbf{G}_l = \mathbb{E}[\tilde{\mathbf{H}}_l(\mathbf{X};\boldsymbol{\Sigma}_\mathbf{X})\mathbf{H}_l^T(\mathbf{X};\boldsymbol{\Sigma}_\mathbf{X})]$ is a $K_{N,l} \times K_{N,l}$ matrix with its $(p,q)$th element

$$G_{l,pq} := \mathbb{E}\left[\tilde{H}_{l,p}(\mathbf{X};\boldsymbol{\Sigma}_\mathbf{X})H_{l,q}(\mathbf{X};\boldsymbol{\Sigma}_\mathbf{X})\right] := \int_{\mathbb{R}^N} \tilde{H}_{l,p}(\mathbf{x};\boldsymbol{\Sigma}_\mathbf{X})H_{l,q}(\mathbf{x};\boldsymbol{\Sigma}_\mathbf{X})\phi_\mathbf{X}(\mathbf{x};\boldsymbol{\Sigma}_\mathbf{X})d\mathbf{x}.$$

From the orthogonality condition (21), $\mathbf{G}_l$ is a diagonal and positive-definite matrix, and hence invertible. Therefore, (22) yields $\mathbf{a}_l = \mathbf{0}$, proving linear independence of the elements of $\mathbf{H}_l(\mathbf{x};\boldsymbol{\Sigma}_\mathbf{X})$ or $\{H_\mathbf{j}(\mathbf{x};\boldsymbol{\Sigma}_\mathbf{X}) : |\mathbf{j}| = l, \mathbf{j} \in \mathbb{N}_0^N\}$. Furthermore, the invertibility of $\mathbf{G}_l$ assures that $\{H_\mathbf{j}(\mathbf{x};\boldsymbol{\Sigma}_\mathbf{X}) : |\mathbf{j}| = l, \mathbf{j} \in \mathbb{N}_0^N\}$ is a spanning set of $\mathcal{V}_l^N$ and, therefore, forms a basis of $\mathcal{V}_l^N$. □

**Corollary 12.** *Let* $\boldsymbol{\Psi}_l(\mathbf{x};\boldsymbol{\Sigma}_\mathbf{X}) := (\Psi_{l,1}(\mathbf{x};\boldsymbol{\Sigma}_\mathbf{X}),\ldots,\Psi_{l,K_{N,l}}(\mathbf{x};\boldsymbol{\Sigma}_\mathbf{X}))^T \in \mathbb{R}^{K_{N,l}}$ *be a column vector constructed from the elements of* $\{\Psi_\mathbf{j}(\mathbf{x};\boldsymbol{\Sigma}_\mathbf{X}) : |\mathbf{j}| = l, \mathbf{j} \in \mathbb{N}_0^N\}$. *Then,* $\mathbf{A}_l = \mathbb{E}[\boldsymbol{\Psi}_l(\mathbf{x};\boldsymbol{\Sigma}_\mathbf{X})\boldsymbol{\Psi}_l^T(\mathbf{x};\boldsymbol{\Sigma}_\mathbf{X})]$, *a* $K_{N,l} \times K_{N,l}$ *matrix with its* $(p,q)$*th element*

$$A_{l,pq} := \mathbb{E}\left[\Psi_{l,p}(\mathbf{X};\boldsymbol{\Sigma}_\mathbf{X})\Psi_{l,q}(\mathbf{X};\boldsymbol{\Sigma}_\mathbf{X})\right] := \int_{\mathbb{R}^N} \Psi_{l,p}(\mathbf{x};\boldsymbol{\Sigma}_\mathbf{X})\Psi_{l,q}(\mathbf{x};\boldsymbol{\Sigma}_\mathbf{X})\phi_\mathbf{X}(\mathbf{x};\boldsymbol{\Sigma}_\mathbf{X})d\mathbf{x},$$

*is symmetric and positive-definite.*

*Proof.* By definition, $\mathbf{A}_l = \mathbf{A}_l^T$. From Proposition 11, the elements of $\boldsymbol{\Psi}_l(\mathbf{x};\boldsymbol{\Sigma}_\mathbf{X})$, a scaled version of $\mathbf{H}_l(\mathbf{x};\boldsymbol{\Sigma}_\mathbf{X})$, are also linearly independent. Therefore, for any $\mathbf{0} \neq \boldsymbol{\alpha}_l \in \mathbb{R}^{K_{N,l}}$, $\boldsymbol{\alpha}_l^T \boldsymbol{\Psi}_l(\mathbf{x};\boldsymbol{\Sigma}_\mathbf{X}) \in \Pi^N$ is a non-zero polynomial, satisfying

$$\boldsymbol{\alpha}_l^T \mathbf{A}_l \boldsymbol{\alpha}_l = \mathbb{E}\left[\left(\boldsymbol{\alpha}_l^T \boldsymbol{\Psi}_l(\mathbf{X};\boldsymbol{\Sigma}_\mathbf{X})\right)^2\right] = \|\boldsymbol{\alpha}_l^T \boldsymbol{\Psi}_l(\mathbf{x};\boldsymbol{\Sigma}_\mathbf{X})\|_{\phi_\mathbf{X} d\mathbf{x}}^2 > 0,$$

as the inner product defined in (5) is positive-definite on $\Pi^N$. Therefore, $\mathbf{A}_l$ is a symmetric, positive-definite matrix. □

**Proposition 13.** *Let* $\mathbf{X} := (X_1,\ldots,X_N)^T : (\Omega,\mathcal{F}) \to (\mathbb{R}^N,\mathcal{B}^N)$, $N \in \mathbb{N}$, *be an* $\mathbb{R}^N$*-valued Gaussian random vector with zero mean; symmetric, positive-definite covariance matrix* $\boldsymbol{\Sigma}_\mathbf{X}$; *and multivariate probability density function* $\phi_\mathbf{X}(\mathbf{x};\boldsymbol{\Sigma}_\mathbf{X})$. *Consistent with the Gaussian measure* $\phi_\mathbf{X}(\mathbf{x})d\mathbf{x}$, *let* $\{H_\mathbf{j}(\mathbf{x};\boldsymbol{\Sigma}_\mathbf{X}) : |\mathbf{j}| = l, \mathbf{j} \in \mathbb{N}_0^N\}$, *the set of multivariate Hermite orthogonal polynomials of degree l, be a basis of* $\mathcal{V}_l^N$. *Then the set of polynomials from the orthogonal sum*

$$\bigoplus_{l \in \mathbb{N}_0} \text{span}\{H_\mathbf{j}(\mathbf{x};\boldsymbol{\Sigma}_\mathbf{X}) : |\mathbf{j}| = l, \mathbf{j} \in \mathbb{N}_0^N\}$$

*is dense in* $L^2(\mathbb{R}^N,\mathcal{B}^N,\phi_\mathbf{X}d\mathbf{x})$. *Moreover,*

$$L^2(\mathbb{R}^N,\mathcal{B}^N,\phi_\mathbf{X}d\mathbf{x}) = \overline{\bigoplus_{l \in \mathbb{N}_0} \mathcal{V}_l^N} \quad (23)$$

*where the overline denotes set closure.*

*Proof.* Define an arbitrary norm $\|\cdot\| : \mathbb{R}^N \to \mathbb{R}_0^+$. According to Skorokhod [22], there exists a real number $\alpha > 0$ such that

$$\int_{\mathbb{R}^N} \exp(\alpha \|\mathbf{x}\|) \phi_{\mathbf{X}}(\mathbf{x}; \boldsymbol{\Sigma}_{\mathbf{X}}) d\mathbf{x} < \infty. \tag{24}$$

In other words, any norm of $\mathbf{x}$ on $\mathbb{R}^N$ is exponentially integrable with respect to the Gaussian probability measure. Now, use Theorem 3.2.18 of Dunkl and Xu [15], which says that if the exponential integrability condition is satisfied, then the space of polynomials $\Pi^N$ is dense in the space $L^2(\mathbb{R}^N, \mathcal{B}^N, \phi_{\mathbf{X}} d\mathbf{x})$. Therefore, the set of polynomial from the orthogonal sum, which is equal to $\Pi^N$ as per (4), is also dense in $L^2(\mathbb{R}^N, \mathcal{B}^N, \phi_{\mathbf{X}} d\mathbf{x})$. Including the limit points of the orthogonal sum yields (23). □

A related subject brought up by Ernst et al. [3] is whether a probability measure is determinate or indeterminate in the Hamburger sense. The multivariate Gaussian probability measure, as it satisfies the exponential integrability condition in (24), is also determinate [15]. In one variable, it is well known that if a measure is determinate, then the space of polynomials is dense in $L^2(\mathbb{R}, \mathcal{B}, \phi_X dx)$. However, this is not universally true for multiple variables. Berg and Thill [23] have shown some rotation-invariant determinate measures for which the spaces of polynomials are not dense. While this matter is not relevant for the Gaussian measure, it can be for non-Gaussian measures.

## 5. Generalized Wiener-Hermite expansion

Let $y(\mathbf{X}) := y(X_1, \ldots, X_N)$ be a real-valued, square-integrable output random variable defined on the same probability space $(\Omega, \mathcal{F}, \mathbb{P})$. The vector space $L^2(\Omega, \mathcal{F}, \mathbb{P})$ is a Hilbert space such that

$$\mathbb{E}\left[y^2(\mathbf{X})\right] := \int_\Omega y^2(\mathbf{X}(\omega)) d\mathbb{P}(\omega) = \int_{\mathbb{R}^N} y^2(\mathbf{x}) \phi_{\mathbf{X}}(\mathbf{x}; \boldsymbol{\Sigma}_{\mathbf{X}}) d\mathbf{x} < \infty$$

with inner product

$$(y(\mathbf{X}), z(\mathbf{X}))_{L^2(\Omega,\mathcal{F},\mathbb{P})} := \int_\Omega y(\mathbf{X}(\omega)) z(\mathbf{X}(\omega)) d\mathbb{P}(\omega) = \int_{\mathbb{R}^N} y(\mathbf{x}) z(\mathbf{x}) \phi_{\mathbf{X}}(\mathbf{x}; \boldsymbol{\Sigma}_{\mathbf{X}}) d\mathbf{x} =: (y(\mathbf{x}), z(\mathbf{x}))_{\phi_{\mathbf{X}} d\mathbf{x}}$$

and norm

$$\|y(\mathbf{X})\|_{L^2(\Omega,\mathcal{F},\mathbb{P})} := \sqrt{(y(\mathbf{X}), y(\mathbf{X}))_{L^2(\Omega,\mathcal{F},\mathbb{P})}} = \sqrt{\mathbb{E}\left[y^2(\mathbf{X})\right]} = \sqrt{(y(\mathbf{x}), y(\mathbf{x}))_{\phi_{\mathbf{X}} d\mathbf{x}}} =: \|y(\mathbf{x})\|_{\phi_{\mathbf{X}} d\mathbf{x}}.$$

It is elementary to show that $y(\mathbf{X}) \in L^2(\Omega, \mathcal{F}, \mathbb{P})$ if and only if $y(\mathbf{x}) \in L^2(\mathbb{R}^N, \mathcal{B}^N, \phi_{\mathbf{X}} d\mathbf{x})$.

### 5.1. Generalized PCE

A generalized PCE of a square-integrable random variable $y(\mathbf{X})$ is simply the expansion of $y(\mathbf{X})$ with respect to an orthogonal polynomial basis of $L^2(\Omega, \mathcal{F}, \mathbb{P})$, formally presented as follows.

**Theorem 14.** *Let $\mathbf{X} := (X_1, \ldots, X_N)^T$, $N \in \mathbb{N}$, be an $\mathbb{R}^N$-valued Gaussian random vector with zero mean, positive-definite covariance matrix $\boldsymbol{\Sigma}_{\mathbf{X}}$, and multivariate probability density function $\phi_{\mathbf{X}}(\mathbf{x}; \boldsymbol{\Sigma}_{\mathbf{X}})$ defined by (1). Then*

*(1) any random variable $y(\mathbf{X}) \in L^2(\Omega, \mathcal{F}, \mathbb{P})$ can be expanded as a Fourier-like infinite series of standardized multivariate Hermite orthogonal polynomials $\{\Psi_{\mathbf{j}}(\mathbf{X}; \boldsymbol{\Sigma}_{\mathbf{X}}) : \mathbf{j} \in \mathbb{N}_0^N\}$, referred to as the generalized PCE of*

$$y(\mathbf{X}) \sim \sum_{\mathbf{j} \in \mathbb{N}_0^N} C_{\mathbf{j}} \Psi_{\mathbf{j}}(\mathbf{X}; \boldsymbol{\Sigma}_{\mathbf{X}}), \tag{25}$$

*where the expansion coefficients $C_{\mathbf{j}} \in \mathbb{R}$, $\mathbf{j} \in \mathbb{N}_0^N$, satisfy the infinite system*

$$\sum_{\substack{\mathbf{k} \in \mathbb{N}_0^N \\ |\mathbf{k}|=|\mathbf{j}|}} C_{\mathbf{k}} \mathbb{E}\left[\Psi_{\mathbf{j}}(\mathbf{X}; \boldsymbol{\Sigma}_{\mathbf{X}}) \Psi_{\mathbf{k}}(\mathbf{X}; \boldsymbol{\Sigma}_{\mathbf{X}})\right] = \mathbb{E}\left[y(\mathbf{X}) \Psi_{\mathbf{j}}(\mathbf{X}; \boldsymbol{\Sigma}_{\mathbf{X}})\right], \, \mathbf{j} \in \mathbb{N}_0^N, \tag{26}$$

*of uncoupled finite-dimensional linear systems; and*

*(2) the generalized PCE of $y(\mathbf{X}) \in L^2(\Omega, \mathcal{F}, \mathbb{P})$ converges to $y(\mathbf{X})$ in mean-square; furthermore, the generalized PCE converges in probability and in distribution.*

*Proof.* If $y(\mathbf{x}) \in L^2(\mathbb{R}^N, \mathcal{B}^N, \phi_\mathbf{X} d\mathbf{x})$, then by Proposition 13, the expansion

$$y(\mathbf{x}) \sim \sum_{l \in \mathbb{N}_0} \text{proj}_l y(\mathbf{x}), \tag{27}$$

with $\text{proj}_l y(\mathbf{x}) : L^2(\mathbb{R}^N, \mathcal{B}^N, \phi_\mathbf{X} d\mathbf{x}) \to \mathcal{V}_l^N$ denoting the projection operator, can be formed. Since standardization is merely scaling, with Proposition 11 in mind, $\mathcal{V}_l^N$ is also spanned by $\{\Psi_\mathbf{j}(\mathbf{x}; \Sigma_\mathbf{X}) : |\mathbf{j}| = l, \mathbf{j} \in \mathbb{N}_0^N\}$. Consequently,

$$\text{proj}_l y(\mathbf{x}) = \sum_{\substack{\mathbf{j} \in \mathbb{N}_0^N \\ |\mathbf{j}|=l}} C_\mathbf{j} \Psi_\mathbf{j}(\mathbf{x}; \Sigma_\mathbf{x}). \tag{28}$$

By definition of the random vector $\mathbf{X}$, the sequence $\{\Psi_\mathbf{j}(\mathbf{X}; \Sigma_\mathbf{x})\}_{\mathbf{j} \in \mathbb{N}_0^N}$ is a basis of $L^2(\Omega, \mathcal{F}, \mathbb{P})$, inheriting the properties of the basis $\{\Psi_\mathbf{j}(\mathbf{x}; \Sigma_\mathbf{x})\}_{\mathbf{j} \in \mathbb{N}_0^N}$ of $L^2(\mathbb{R}^N, \mathcal{B}^N, \phi_\mathbf{X} d\mathbf{x})$. Therefore, (27) and (28) lead to the expansion in (25).

In reference to Proposition 13, recognize that the set of polynomials from the orthogonal sum

$$\bigoplus_{l \in \mathbb{N}_0} \text{span}\{\Psi_\mathbf{j}(\mathbf{x}; \Sigma_\mathbf{X}) : |\mathbf{j}| = l, \mathbf{j} \in \mathbb{N}_0^N\} = \Pi^N \tag{29}$$

is also dense in $L^2(\mathbb{R}^N, \mathcal{B}^N, \phi_\mathbf{X} d\mathbf{x})$. Therefore, one has the Bessel's inequality [24]

$$\mathbb{E}\left[\sum_{\mathbf{j} \in \mathbb{N}_0^N} C_\mathbf{j} \Psi_\mathbf{j}(\mathbf{X}; \Sigma_\mathbf{X})\right]^2 \leq \mathbb{E}\left[y^2(\mathbf{X})\right],$$

proving that the generalized PCE converges in mean-square or $L^2$. To determine the limit of convergence, invoke again Proposition 13, which implies that the set $\{\Psi_\mathbf{j}(\mathbf{x}; \Sigma_\mathbf{X}) : \mathbf{j} \in \mathbb{N}_0^N\}$ is complete in $L^2(\mathbb{R}^N, \mathcal{B}^N, \phi_\mathbf{X} d\mathbf{x})$. Therefore, Bessel's inequality becomes an equality

$$\mathbb{E}\left[\sum_{\mathbf{j} \in \mathbb{N}_0^N} C_\mathbf{j} \Psi_\mathbf{j}(\mathbf{X}; \Sigma_\mathbf{X})\right]^2 = \mathbb{E}\left[y^2(\mathbf{X})\right],$$

known as the Parseval identity [24] for a multivariate orthogonal system, for every random variable $y(\mathbf{X}) \in L^2(\Omega, \mathcal{F}, \mathbb{P})$. Furthermore, as the PCE converges in mean-square, it does so in probability. Moreover, as the expansion converges in probability, it also converges in distribution.

Finally, to find the expansion coefficients, define a second moment

$$e_{\text{PCE}} := \mathbb{E}\left[y(\mathbf{X}) - \sum_{\mathbf{k} \in \mathbb{N}_0^N} C_\mathbf{k} \Psi_\mathbf{k}(\mathbf{X}; \Sigma_\mathbf{X})\right]^2 \tag{30}$$

of the difference between $y(\mathbf{X})$ and its full PCE. Differentiate both sides of (30) with respect to $C_\mathbf{j}$, $\mathbf{j} \in \mathbb{N}_0^N$, to write

$$
\begin{aligned}
\frac{\partial e_{\text{PCE}}}{\partial C_\mathbf{j}} &= \frac{\partial}{\partial C_\mathbf{j}} \mathbb{E}\left[ y(\mathbf{X}) - \sum_{\mathbf{k} \in \mathbb{N}_0^N} C_\mathbf{k} \Psi_\mathbf{k}(\mathbf{X}; \Sigma_\mathbf{X}) \right]^2 \\
&= \mathbb{E}\left[ \frac{\partial}{\partial C_\mathbf{j}} \left\{ y(\mathbf{X}) - \sum_{\mathbf{k} \in \mathbb{N}_0^N} C_\mathbf{k} \Psi_\mathbf{k}(\mathbf{X}; \Sigma_\mathbf{X}) \right\}^2 \right] \\
&= 2\mathbb{E}\left[ \left\{ \sum_{\mathbf{k} \in \mathbb{N}_0^N} C_\mathbf{k} \Psi_\mathbf{k}(\mathbf{X}; \Sigma_\mathbf{X}) - y(\mathbf{X}) \right\} \Psi_\mathbf{j}(\mathbf{X}; \Sigma_\mathbf{X}) \right] \\
&= 2\left\{ \sum_{\mathbf{k} \in \mathbb{N}_0^N} C_\mathbf{k} \mathbb{E}\left[ \Psi_\mathbf{j}(\mathbf{X}; \Sigma_\mathbf{X}) \Psi_\mathbf{k}(\mathbf{X}; \Sigma_\mathbf{X}) \right] - \mathbb{E}\left[ y(\mathbf{X}) \Psi_\mathbf{j}(\mathbf{X}; \Sigma_\mathbf{X}) \right] \right\} \\
&= 2\left\{ \sum_{\substack{\mathbf{k} \in \mathbb{N}_0^N \\ |\mathbf{k}|=|\mathbf{j}|}} C_\mathbf{k} \mathbb{E}\left[ \Psi_\mathbf{j}(\mathbf{X}; \Sigma_\mathbf{X}) \Psi_\mathbf{k}(\mathbf{X}; \Sigma_\mathbf{X}) \right] - \mathbb{E}\left[ y(\mathbf{X}) \Psi_\mathbf{j}(\mathbf{X}; \Sigma_\mathbf{X}) \right] \right\}.
\end{aligned}
\quad (31)
$$

Here, the second, third, fourth, and last lines are obtained by interchanging the differential and expectation operators, performing the differentiation, swapping the expectation and summation operators, and applying Corollary 9, respectively. The interchanges are permissible as the infinite sum is convergent as demonstrated in the preceding paragraph. Setting $\partial e_{\text{PCE}}/\partial C_\mathbf{j} = 0$ in (31) yields (26), completing the proof. □

The linear system (26) can also be derived by simply replacing $y(\mathbf{X})$ in (26) with the full PCE and then using Corollary 9. In contrast, the proof given here demonstrates that the PCE coefficients are determined optimally.

The generalized PCE presented here should not be confused with that of Xiu and Karniadakis [5]. The generalization in this work extends the applicability of the classical Wiener-Hermite PCE for arbitrary but dependent Gaussian probability distributions of random input. In contrast, the existing generalized PCE [5] still requires independence of random input, but can account for non-Gaussian variables, provided that the marginal probability measures are determinate.

**Corollary 15.** *If $\mathbf{X} = (X_1, \ldots, X_N)^T$ comprises independent Gaussian random variables, each with zero mean and unit variance, then $\Sigma_\mathbf{X} = \mathbf{I}$ and $\Psi_\mathbf{j}(\mathbf{x}; \Sigma_\mathbf{X}) = \Pi_{i=1}^N \Psi_{j_i}(x_i)$ with $\Psi_{j_i}(x_i)$ representing the $j_i$th-degree univariate Hermite orthonormal polynomial in $x_i$. In which case, the generalized PCE reduces to the classical PCE, yielding*

$$y(\mathbf{X}) \sim \sum_{\mathbf{j} \in \mathbb{N}_0^N} C_\mathbf{j} \prod_{i=1}^N \Psi_{j_i}(X_i)$$

*with the expansion coefficients*

$$C_\mathbf{j} = \mathbb{E}\left[ y(\mathbf{X}) \prod_{i=1}^N \Psi_{j_i}(X_i) \right]. \quad (32)$$

Note that the linear system (26) of the generalized PCE is coupled with respect to the coefficients of the same degree. This is due to weak orthogonality of Hermite polynomials for dependent variables. The Hermite polynomials for independent variables, by contrast, are strongly orthogonal. In consequence, there are no such interactions among respective coefficients, as presented in (32), for the classical PCE.

It should be emphasized that the function $y$ must be square-integrable for the mean-square and other convergences to hold. However, the rate of convergence depends on the smoothness of the function. The smoother the function, the faster the convergence. If the function is a polynomial, then its generalized PCE exactly reproduces the function. These well-known results from the literature of classical PCE extend to the generalized PCE and can be proved using classical approximation theory.

Note that the infinite series in (25) does not necessarily converge almost surely to $y(\mathbf{X})$, that is, for $m \in \mathbb{N}_0$, $\sum_{\mathbf{j} \in \mathbb{N}_0^N, |\mathbf{j}| \leq m} C_\mathbf{j} \Psi_\mathbf{j}(\mathbf{X}(\omega); \Sigma_\mathbf{X})$ may not approach $y(\mathbf{X}(\omega))$ as $m \to \infty$. Furthermore, it is not guaranteed that the moments of PCE of order larger than two will converge. These known fundamental limitations of classical PCE persist in the generalized PCE.

### 5.2. Truncation

The generalized PCE contains an infinite number of orthogonal polynomials or coefficients. In practice, the number must be finite, meaning that the PCE must be truncated. But there are multiple ways to perform the truncation. A popular approach, adopted in this work, entails retaining all polynomials with the total degree $|\mathbf{j}|$ less than or equal to $m \in \mathbb{N}$. The result is an $m$th-order generalized PCE approximation [3]

$$y_m(\mathbf{X}) = \sum_{\substack{\mathbf{j} \in \mathbb{N}_0^N \\ 0 \le |\mathbf{j}| \le m}} C_{\mathbf{j}} \Psi_{\mathbf{j}}(\mathbf{X}; \boldsymbol{\Sigma}_{\mathbf{X}}) = \sum_{l=0}^{m} \sum_{\substack{\mathbf{j} \in \mathbb{N}_0^N \\ |\mathbf{j}|=l}} C_{\mathbf{j}} \Psi_{\mathbf{j}}(\mathbf{X}; \boldsymbol{\Sigma}_{\mathbf{X}}) \tag{33}$$

of $y(\mathbf{X})$, which contains

$$L_{N,m} = \binom{N+m}{m} = \frac{(N+m)!}{N!m!}$$

number of expansion coefficients, satisfying the finite-dimensional linear system

$$\sum_{\substack{\mathbf{k} \in \mathbb{N}_0^N \\ |\mathbf{k}|=|\mathbf{j}|}} C_{\mathbf{k}} \mathbb{E}\left[\Psi_{\mathbf{j}}(\mathbf{X}; \boldsymbol{\Sigma}_{\mathbf{X}}) \Psi_{\mathbf{k}}(\mathbf{X}; \boldsymbol{\Sigma}_{\mathbf{X}})\right] = \mathbb{E}\left[y(\mathbf{X}) \Psi_{\mathbf{j}}(\mathbf{X}; \boldsymbol{\Sigma}_{\mathbf{X}})\right], \quad 0 \le |\mathbf{j}| \le m. \tag{34}$$

It is natural to ask about the approximation quality of (33). Since the set $\{\Psi_{\mathbf{j}}(\mathbf{x}; \boldsymbol{\Sigma}_{\mathbf{X}}) : \mathbf{j} \in \mathbb{N}_0^N\}$ or $\{\Psi_{\mathbf{j}}(\mathbf{X}; \boldsymbol{\Sigma}_{\mathbf{X}}) : \mathbf{j} \in \mathbb{N}_0^N\}$ is complete in $L^2(\mathbb{R}^N, \mathcal{B}^N, \phi_{\mathbf{X}} d\mathbf{x})$ or $L^2(\Omega, \mathcal{F}, \mathbb{P})$, the truncation error $y(\mathbf{X}) - y_m(\mathbf{X})$ is orthogonal to any element of the space from which $y_m(\mathbf{X})$ is chosen, as demonstrated below.

**Corollary 16.** *The truncation error $y(\mathbf{X}) - y_m(\mathbf{X})$ is orthogonal to the span of $\{\Psi_{\mathbf{j}}(\mathbf{X}; \boldsymbol{\Sigma}_{\mathbf{X}}), 0 \le |\mathbf{j}| \le m\}$. Moreover, $\mathbb{E}[y(\mathbf{X}) - y_m(\mathbf{X})]^2 \to 0$ as $m \to \infty$.*

*Proof.* Let

$$\bar{y}_m(\mathbf{X}) := \sum_{\substack{\mathbf{k} \in \mathbb{N}_0^N \\ 0 \le |\mathbf{k}| \le m}} \bar{C}_{\mathbf{k}} \Psi_{\mathbf{k}}(\mathbf{X}; \boldsymbol{\Sigma}_{\mathbf{X}}),$$

with arbitrary expansion coefficients $\bar{C}_{\mathbf{k}}$, $0 \le |\mathbf{k}| \le m$, be any element of the subspace of $L^2(\Omega, \mathcal{F}, \mathbb{P})$ spanned by $\{\Psi_{\mathbf{k}}(\mathbf{X}; \boldsymbol{\Sigma}_{\mathbf{X}}) : 0 \le |\mathbf{k}| \le m\}$. Then

$$\begin{aligned}
&= \mathbb{E}\left[\{y(\mathbf{X}) - y_m(\mathbf{X})\} \bar{y}_m(\mathbf{X})\right] \\
&= \mathbb{E}\left[\left\{\sum_{\substack{\mathbf{j} \in \mathbb{N}_0^N \\ m+1 \le |\mathbf{j}| < \infty}} C_{\mathbf{j}} \Psi_{\mathbf{j}}(\mathbf{X}; \boldsymbol{\Sigma}_{\mathbf{X}})\right\} \left\{\sum_{\substack{\mathbf{k} \in \mathbb{N}_0^N \\ 0 \le |\mathbf{k}| \le m}} \bar{C}_{\mathbf{k}} \Psi_{\mathbf{k}}(\mathbf{X}; \boldsymbol{\Sigma}_{\mathbf{X}})\right\}\right] \\
&= \sum_{\substack{\mathbf{j},\mathbf{k} \in \mathbb{N}_0^N \\ m+1 \le |\mathbf{j}| < \infty \\ 0 \le |\mathbf{k}| \le m}} C_{\mathbf{j}} \bar{C}_{\mathbf{k}} \mathbb{E}\left[\Psi_{\mathbf{j}}(\mathbf{X}; \boldsymbol{\Sigma}_{\mathbf{X}}) \Psi_{\mathbf{k}}(\mathbf{X}; \boldsymbol{\Sigma}_{\mathbf{X}})\right] \\
&= 0,
\end{aligned}$$

where the last line follows from Corollary 9, proving the first part of the proposition. For the latter part, the Pythagoras theorem yields

$$\mathbb{E}[\{y(\mathbf{X}) - y_m(\mathbf{X})\}^2] + \mathbb{E}[y_m^2(\mathbf{X})] = \mathbb{E}[y(\mathbf{X})^2].$$

From Theorem 14, $\mathbb{E}[y_m^2(\mathbf{X})] \to \mathbb{E}[y^2(\mathbf{X})]$ as $m \to \infty$. Therefore, $\mathbb{E}[\{y(\mathbf{X}) - y_m(\mathbf{X})\}^2] \to 0$ as $m \to \infty$. □

---

[3]The nouns *degree* and *order* associated with PCE or Hermite polynomials are used synonymously in the paper.

The second part of Corollary 16 entails $L^2$ convergence, which is the same as the mean-square convergence described in Theorem 14. However, an alternative route is chosen for the proof of Corollary 16.

*5.2.1. Second-moment statistics*

The $m$th-order generalized PCE approximation $y_m(\mathbf{X})$ can be viewed as a surrogate of $y(\mathbf{X})$. Therefore, relevant probabilistic characteristics of $y(\mathbf{X})$, including its first two moments and probability density function, if it exists, can be estimated from the statistical properties of $y_m(\mathbf{X})$.

Applying the expectation operator on $y_m(\mathbf{X})$ and $y(\mathbf{X})$ in (25) and (33) and imposing Corollary 9, their means

$$\mathbb{E}[y_m(\mathbf{X})] = \mathbb{E}[y(\mathbf{X})] = C_{\mathbf{0}} \tag{35}$$

are the same as the *zero*-degree expansion coefficient and are independent of $m$. Therefore, the generalized PCE truncated for any value of $m$ yields the exact mean. The formulae for the means in the classical and generalized PCE are the same, although the respective expansion coefficients involved are not. Nonetheless, $\mathbb{E}[y_m(\mathbf{X})]$ will be referred to as the $m$th-order generalized PCE approximation of the mean of $y(\mathbf{X})$.

Applying the expectation operator again, this time on $[y_m(\mathbf{X}) - C_{\mathbf{0}}]^2$ and $[y(\mathbf{X}) - C_{\mathbf{0}}]^2$, and employing Corollary 9 results in the variances

$$\text{var}[y_m(\mathbf{X})] = \sum_{\substack{\mathbf{j} \in \mathbb{N}_0^N \\ 1 \le |\mathbf{j}| \le m}} C_{\mathbf{j}}^2 + \sum_{\substack{\mathbf{j},\mathbf{k} \in \mathbb{N}_0^N \\ 1 \le |\mathbf{j}|,|\mathbf{k}| \le m \\ |\mathbf{j}|=|\mathbf{k}|, \mathbf{j} \ne \mathbf{k}}} C_{\mathbf{j}} C_{\mathbf{k}} \mathbb{E}\left[\Psi_{\mathbf{j}}(\mathbf{X}; \Sigma_{\mathbf{X}}) \Psi_{\mathbf{k}}(\mathbf{X}; \Sigma_{\mathbf{X}})\right] \tag{36}$$

and

$$\text{var}[y(\mathbf{X})] = \sum_{\substack{\mathbf{j} \in \mathbb{N}_0^N \\ 1 \le |\mathbf{j}| < \infty}} C_{\mathbf{j}}^2 + \sum_{\substack{\mathbf{j},\mathbf{k} \in \mathbb{N}_0^N \\ 1 \le |\mathbf{j}|,|\mathbf{k}| < \infty \\ |\mathbf{j}|=|\mathbf{k}|, \mathbf{j} \ne \mathbf{k}}} C_{\mathbf{j}} C_{\mathbf{k}} \mathbb{E}\left[\Psi_{\mathbf{j}}(\mathbf{X}; \Sigma_{\mathbf{X}}) \Psi_{\mathbf{k}}(\mathbf{X}; \Sigma_{\mathbf{X}})\right]$$

of $y_m(\mathbf{X})$ and $y(\mathbf{X})$, respectively. The condition $1 \le |\mathbf{j}|, |\mathbf{k}| \le m$ in the summation means $1 \le |\mathbf{j}| \le m$ and $1 \le |\mathbf{k}| \le m$. In (36), the lower limit of $|\mathbf{j}|$ exceeds the upper limit when $m = 0$, yielding $\text{var}[y_0(\mathbf{X})] = 0$. This is consistent with $y_0(\mathbf{X}) = C_{\mathbf{0}}$, a constant function producing no variance. Clearly, $\text{var}[y_m(\mathbf{X})]$, referred to as the $m$th-order generalized PCE approximation of the variance of $y(\mathbf{X})$, approaches $\text{var}[y(\mathbf{X})]$, the exact variance of $y(\mathbf{X})$, as $m \to \infty$. Compared with the classical PCE, the formulae for the variances in the generalized PCE include a second sum, which represents the contribution from the correlation properties of input variables $\mathbf{X}$. The second sum vanishes in the formulae for the variances in the classical PCE as $\mathbf{X}$ comprises only independent variables.

Being convergent in probability and distribution, the probability density function of $y(\mathbf{X})$, if it exists, can also be estimated by that of $y_m(\mathbf{X})$. However, no analytical formula exists for the density function. In that case, the density can be estimated by MCS of $y_m(\mathbf{X})$. Such simulation should not be confused with crude MCS of $y(\mathbf{X})$, commonly used for producing benchmark results whenever possible. The crude MCS can be expensive or even prohibitive, particularly when the sample size needs to be very large for estimating tail probabilistic characteristics. In contrast, the MCS embedded in the generalized PCE approximation requires evaluations of simple polynomial functions that describe $y_m$. Therefore, a relatively large sample size can be accommodated in the PCE approximation even when $y$ is expensive to evaluate.

*5.2.2. Expansion coefficients*

According to (34), determining the expansion coefficients of the $m$th-order generalized PCE approximation requires solving an $(L_{N,m} \times L_{N,m})$ system of linear equations. However, the coefficients interact with each other only for a specific degree. Therefore, the coefficients for each degree can be determined independently, described as follows.

Let $0 \le l \le m$ be a degree of orthogonal polynomials for which there are

$$K_{N,l} = \binom{N+l-1}{l} = \frac{(N+l-1)!}{l!(N-1)!}$$

$l$th-degree expansion coefficients $C_\mathbf{j}$, $|\mathbf{j}| = l$. To determine all $l$th-degree coefficients, only a $(K_{N,l} \times K_{N,l})$ linear system,

$$\sum_{\substack{\mathbf{k} \in \mathbb{N}_0^N \\ |\mathbf{k}|=|\mathbf{j}|}} C_\mathbf{k} \mathbb{E}\left[\Psi_\mathbf{j}(\mathbf{X}; \Sigma_\mathbf{X})\Psi_\mathbf{k}(\mathbf{X}; \Sigma_\mathbf{X})\right] = \mathbb{E}\left[y(\mathbf{X})\Psi_\mathbf{j}(\mathbf{X}; \Sigma_\mathbf{X})\right], \quad |\mathbf{j}| = l, \tag{37}$$

has to be solved. Appendix A gives further details on how to build the matrix form of the linear system. When (37) is solved for $l = 0, \ldots, m$, then all $L_{N,m}$ expansion coefficients for degree at most $m$ have been determined. Obviously, $L_{N,m} = \sum_{l=0}^{m} K_{N,l}$.

The linear system (37) requires calculating the expectations $\mathbb{E}[y(\mathbf{X})\Psi_\mathbf{j}(\mathbf{X}; \Sigma_\mathbf{X})]$ for $|\mathbf{j}| = l$. These expectations are various $N$-dimensional integrals on $\mathbb{R}^N$, which cannot be determined analytically or exactly if $y$ is a general function. Furthermore, for large $N$, a full numerical integration employing an $N$-dimensional tensor product of a univariate quadrature formula is computationally expensive and likely prohibitive. Therefore, alternative means of estimating these expectations or integrals must be pursued. One approach entails exploiting smart combinations of low-dimensional numerical integrations, such as sparse-grid quadrature [25] and dimension-reduction integration [26], to approximate a high-dimensional integral. The other approach consists of efficient sampling methods, such as quasi Monte Carlo simulation (QMCS) [27], importance sampling with Monte Carlo [28], and Markov chain Monte Carlo [29], to name a few. In the latter approach, one hopes to attain sufficiently accurate estimates of the expansion coefficients for a relatively low sample size. However, if the sample size required is too high, then the statistics of $y(\mathbf{X})$ can be estimated directly, raising a question about the need for a PCE approximation in the first place. The topic merits further study.

### 5.2.3. Numerical implementation

Algorithm 1 describes a procedure for developing an $m$th-order generalized PCE approximation $y_m(\mathbf{X})$ of a general square-integrable function $y(\mathbf{X})$. It includes calculation of the mean and variance of $y_m(\mathbf{X})$.

---

**Algorithm 1:** Generalized PCE approximation and second-moment statistics

**Input**: The total number $N$ of Gaussian input variables $\mathbf{X} = (X_1, \ldots, X_N)^T$, a positive-definite covariance matrix $\Sigma_\mathbf{X}$ of $\mathbf{X}$, a square-integrable function $y(\mathbf{X})$, and the largest order $m$ of orthogonal polynomials

**Output**: The $m$th-order PCE approximation $y_m(\mathbf{X})$ of $y(\mathbf{X})$, mean and variance of $y_m(\mathbf{X})$

1. **for** $l \leftarrow 0$ *to* $m$ **do**
2.     Generate Hermite polynomials $H_\mathbf{j}(\mathbf{x}; \Sigma_\mathbf{X})$ and $\Psi_\mathbf{j}(\mathbf{x}; \Sigma_\mathbf{X})$, $|\mathbf{j}| = l$
       `/* from (6) and (7) */`
3.     Calculate $\mathbb{E}[H_\mathbf{j}(\mathbf{X}; \Sigma_\mathbf{X})H_\mathbf{k}(\mathbf{X}; \Sigma_\mathbf{X})]$ and $\mathbb{E}[H_\mathbf{j}^2(\mathbf{X}; \Sigma_\mathbf{X})]$, $|\mathbf{j}| = |\mathbf{k}| = l$
       `/* from (12) and (13) */`
4.     Calculate $\mathbb{E}[\Psi_\mathbf{j}(\mathbf{X}; \Sigma_\mathbf{X})\Psi_\mathbf{k}(\mathbf{X}; \Sigma_\mathbf{X})]$, $|\mathbf{j}| = |\mathbf{k}| = l$
       `/* from (19) */`
5.     Calculate or estimate $\mathbb{E}[y(\mathbf{X})\Psi_\mathbf{j}(\mathbf{X}; \Sigma_\mathbf{X})]$, $|\mathbf{j}| = l$
       `/* from reduced integration or sampling methods */`
6.     Construct the system matrix $\mathbf{A}_l$ and vector $\mathbf{b}_l$
       `/* from Appendix A */`
7.     Solve the linear system $\mathbf{A}_l \mathbf{c}_l = \mathbf{b}_l$ for $l$th-order PCE coefficients
       `/* from Appendix A */`
8. Compile a set $\{C_\mathbf{j}, 0 \leq |\mathbf{j}| \leq m\}$ of at most $m$th-order PCE coefficients and hence construct the $m$th-order PCE approximation $y_m(\mathbf{X})$
       `/* from (33) */`
9. Calculate the mean $\mathbb{E}[y_m(\mathbf{X})]$ and variance $\text{var}[y_m(\mathbf{X})]$
       `/* from (35) and (36) */`

---

When the covariance matrix is positive-definite, as assumed here, the Cholesky factorization of the covariance matrix leads to a linear map between dependent and independent Gaussian variables. Therefore, the classical PCE can also be used for tackling dependent Gaussian variables. In contrast, the generalized PCE proposed provides an alternative means of solving stochastic problems with dependent Gaussian variables directly, that is, without the transformation. More importantly, if the input variables are both dependent and non-Gaussian, then the Cholesky factorization is inadequate, if not useless, and the map becomes nonlinear in general, rendering the classical PCE inefficient. In which case, the use of multivariate orthogonal polynomials and generalized PCE, if they exist, is more relevant and perhaps necessary. The extension to non-Gaussian variables is discussed in the last subsection.

*5.3. Infinitely many input variables*

In uncertainty quantification, information theory, and stochastic process, functions depending on a countable sequence $\{X_i\}_{i\in\mathbb{N}}$ of input random variables need to be considered. Does the generalized PCE proposed still apply as in the case of finitely many random variables? The following proposition provides the answer.

**Proposition 17.** *Let $\{X_i\}_{i\in\mathbb{N}}$ be a countable sequence of Gaussian random variables defined on the probability space $(\Omega, \mathcal{F}_\infty, \mathbb{P})$, where $\mathcal{F}_\infty := \sigma(\{X_i\}_{i\in\mathbb{N}})$ is the associated $\sigma$-algebra generated. Then the generalized PCE of $y(\{X_i\}_{i\in\mathbb{N}}) \in L^2(\Omega, \mathcal{F}_\infty, \mathbb{P})$, where $y : \mathbb{R}^\mathbb{N} \to \mathbb{R}$, converges to $y(\{X_i\}_{i\in\mathbb{N}})$ in mean-square. Moreover, the generalized PCE converges in probability and in distribution.*

*Proof.* According to Proposition 13, $\Pi^N$ is dense in $L^2(\mathbb{R}^N, \mathcal{B}^N, \phi_\mathbf{X} d\mathbf{x})$ and hence in $L^2(\Omega, \mathcal{F}_N, \mathbb{P})$ for every $N \in \mathbb{N}$, where $\mathcal{F}_N := \sigma(\{X_i\}_{i=1}^N)$ is the associated $\sigma$-algebra generated by $\{X_i\}_{i=1}^N$.[4] Now, apply Theorem 3.8 of Ernst et al. [3], which says that if $\Pi^N$ is dense in $L^2(\Omega, \mathcal{F}_N, \mathbb{P})$ for every $N \in \mathbb{N}$, then

$$\Pi^\infty := \bigcup_{N=1}^\infty \Pi^N,$$

a subspace of $L^2(\Omega, \mathcal{F}_\infty, \mathbb{P})$, is also dense in $L^2(\Omega, \mathcal{F}_\infty, \mathbb{P})$. But, using (29),

$$\Pi^\infty = \bigcup_{N=1}^\infty \bigoplus_{l\in\mathbb{N}_0} \text{span}\{\Psi_\mathbf{j} : |\mathbf{j}| = l, \mathbf{j} \in \mathbb{N}_0^N\} = \bigcup_{N=1}^\infty \text{span}\{\Psi_\mathbf{j} : \mathbf{j} \in \mathbb{N}_0^N\},$$

demonstrating that the set of polynomials from the union is dense in $L^2(\Omega, \mathcal{F}_\infty, \mathbb{P})$. Therefore, the generalized PCE of $y(\{X_i\}_{i\in\mathbb{N}}) \in L^2(\Omega, \mathcal{F}_\infty, \mathbb{P})$ converges to $y(\{X_i\}_{i\in\mathbb{N}})$ in mean-square. Since the mean-square convergence is stronger than the convergence in probability or in distribution, the latter modes of convergence follow readily. □

*5.4. Extension for non-Gaussian measures*

Although the paper focuses on PCE for Gaussian measures, a further generalization is possible for non-Gaussian measures. However, a few important conditions must be fulfilled before proceeding with the generalization. First and foremost, the non-Gaussian measures must be determinate. More importantly, the set of orthogonal polynomials consistent with a non-Gaussian measure, if they exist, must be dense or complete in $L^2(\Omega, \mathcal{F}, \mathbb{P})$. Otherwise, the resultant PCE may not converge to the correct limit. It is important to note that the denseness condition is easily satisfied for a probability density function with a compact support. For an unbounded support, the exponential integrability of a norm, as done here for the Gaussian density function, or other alternatives will have to be established.

Second, numerical methods must be used in general to generate measure-consistent orthogonal polynomials. In this case, the Gram-Schmidt orthogonalization [12], commonly used for building univariate polynomials, is useful for constructing multivariate polynomials as well. However, an important difference between univariate polynomials and multivariate polynomials is the lack of an obvious natural order in the latter. The natural order for monomials of univariate polynomials is the degree order; that is, one orders monomials according to their degree. For multivariate polynomials, there are many options, such as lexicographic order, graded lexicographic order, and reversed graded lexicographic order, to name just three. There is no natural choice, and different orders will give different sequences of orthogonal polynomials from the Gram-Schmidt orthogonalization.

---

[4] With a certain abuse of notation, $\Pi^N$ is used here as a set of polynomial functions of both real variables ($\mathbf{x}$) and random variables ($\mathbf{X}$).

Last but not least, deriving an analytical formula for the second-moment properties of orthogonal polynomials for arbitrary non-Gaussian measures is nearly impossible. Having said so, these properties, which represent high-dimensional integrals comprising products of orthogonal polynomials, can be estimated by numerical integration with an arbitrary precision even when $N$ is large. This is because no generally expensive output function evaluations are involved. Given that these issues are properly accounted for, the rest of the PCE proposed should work for non-Gaussian measures.

## 6. Numerical examples

Three examples, involving an explicit polynomial function, an implicit non-polynomial function satisfying a stochastic ordinary differential equation, and an implicit function derived from finite-element random eigenvalue analysis, are presented to illustrate the generalized PCE.

### 6.1. Example 1

As introduced in the author's earlier work [30], consider a symmetric, quadratic, polynomial function

$$y(\mathbf{X}) = 12 + 4X_1 + 4X_2 + 4X_3 + X_1X_2 + X_1X_3 + X_2X_3$$

of a trivariate Gaussian random vector $\mathbf{X} = (X_1, X_2, X_3)^T$, which has mean $\boldsymbol{\mu}_\mathbf{X} = \mathbb{E}[\mathbf{X}] = \mathbf{0} \in \mathbb{R}^3$, positive-definite covariance matrix

$$\boldsymbol{\Sigma}_\mathbf{X} = \mathbb{E}\left[\mathbf{X}\mathbf{X}^T\right] = \begin{bmatrix} \sigma_1^2 & \rho_{12}\sigma_1\sigma_2 & \rho_{13}\sigma_1\sigma_3 \\ & \sigma_2^2 & \rho_{23}\sigma_2\sigma_3 \\ \text{(sym.)} & & \sigma_3^2 \end{bmatrix} \in \mathbb{S}_+^3,$$

comprising variances $\sigma_i^2 = 1$ of $X_i$ for $i = 1, 2, 3$ and correlation coefficients $\rho_{ij}$ between $X_i$ and $X_j$, $i, j = 1, 2, 3$, $i \neq j$, and a joint probability density function described by (1) for $N = 3$. Four cases of correlation coefficients with varied strengths and types of statistical dependence among random variables were examined: (1) $\rho_{12} = \rho_{13} = \rho_{23} = 0$ (no correlation); (2) $\rho_{12} = \rho_{13} = \rho_{23} = 1/5$ (equal correlation); (3) $\rho_{12} = 1/5$, $\rho_{13} = 2/5$, $\rho_{23} = 4/5$ (positive correlation); and (4) $\rho_{12} = -1/5$, $\rho_{13} = 2/5$, $\rho_{23} = -4/5$ (mixed correlation). The objective of this example is to explain the construction of the generalized PCE and the calculation of the second-moment statistics for all four cases of correlation coefficients.

Since $y(\mathbf{X})$ is a quadratic polynomial, its second-order generalized PCE approximation

$$y_2(\mathbf{X}) = \sum_{\substack{\mathbf{j} \in \mathbb{N}_0^3 \\ 0 \leq |\mathbf{j}| \leq 2}} C_\mathbf{j} \Psi_\mathbf{j}(\mathbf{X}; \boldsymbol{\Sigma}_\mathbf{X}) = \sum_{l=0}^{2} \sum_{\substack{\mathbf{j} \in \mathbb{N}_0^3 \\ |\mathbf{j}|=l}} C_\mathbf{j} \Psi_\mathbf{j}(\mathbf{X}; \boldsymbol{\Sigma}_\mathbf{X})$$

was built following Algorithm 1 to reproduce the former. Given $N = 3$ and $m = 2$, the number of multivariate Hermite polynomials or PCE coefficients is $L_{3,2} = (3 + 2)!/(3!2!) = 10$. Table 1 presents all ten standardized Hermite orthogonal polynomials, obtained using (6) and (7), for four distinct cases of correlation coefficients. The corresponding expansion coefficients were calculated by forming the system matrix $\mathbf{A}_l$ and vector $\mathbf{b}_l$, as explained in Appendix A, and then solving the linear system for the vector $\mathbf{c}_l$ of coefficients for each degree $l = 0, 1, 2$ separately. While the expectations involved in $\mathbf{A}_l$ were determined from the proposed analytical formulae described by (19) and (20), the expectations contained in $\mathbf{b}_l$ were obtained by analytical integrations, which is possible for the function $y$ chosen. Therefore, all coefficients of the generalized PCE, listed in Table 2, were determined exactly.

From Tables 1 and 2, clearly, the orthogonal polynomials and expansion coefficients vary with the correlation structure, but when added together they reconstruct the same function $y$ whether or not the random variables are independent. When there is no correlation between any two random variables (Case 1), the standardized multivariate orthogonal polynomials are products of univariate orthonormal polynomials, and the expansion coefficients are merely the coefficients of the function $y$, as expected in the classical PCE. In other words, the generalized PCE reduces to the classical PCE for independent random variables. For the remaining three cases (Cases 2 through 4), the zeroth-order orthogonal polynomials are equal to *one*, the same constant as in Case 1, but the corresponding expansion coefficients vary with the correlation properties. Moreover, the first- and second-order orthogonal polynomials contain additional terms of the same degree that are not present in the original function to begin with. It is easy to verify from Corollaries 9 and 10 that all first- and second-order orthogonal polynomials have *zero* means and are either strongly orthogonal for Case 1 or weakly orthogonal for Cases 2 through 4.

Table 1: Zeroth-, first-, and second-order standardized Hermite orthogonal polynomials in Example 1.

| Case 1: $\rho_{12} = \rho_{13} = \rho_{23} = 0$ | Case 2: $\rho_{12} = \rho_{13} = \rho_{23} = 1/5$ |
|---|---|
| $\Psi_{(0,0,0)} = 1$ | $\Psi_{(0,0,0)} = 1$ |
| $\Psi_{(1,0,0)} = x_1$ | $\Psi_{(1,0,0)} = \frac{1}{2}\sqrt{\frac{5}{42}}(6x_1 - x_2 - x_3)$ |
| $\Psi_{(0,1,0)} = x_2$ | $\Psi_{(0,1,0)} = -\frac{1}{2}\sqrt{\frac{5}{42}}(x_1 - 6x_2 + x_3)$ |
| $\Psi_{(0,0,1)} = x_3$ | $\Psi_{(0,0,1)} = -\frac{1}{2}\sqrt{\frac{5}{42}}(x_1 + x_2 - 6x_3)$ |
| $\Psi_{(2,0,0)} = \frac{1}{\sqrt{2}}(x_1^2 - 1)$ | $\Psi_{(2,0,0)} = [180x_1^2 - 60x_1(x_2 + x_3) + 5x_2^2 + 5x_3^2 + 10x_2x_3 - 168]/168\sqrt{2}$ |
| $\Psi_{(1,1,0)} = x_1 x_2$ | $\Psi_{(1,1,0)} = -[30x_1^2 - 5x_1(37x_2 - 5x_3) + 30x_2^2 - 5x_3^2 + 25x_2x_3 - 28]/28\sqrt{37}$ |
| $\Psi_{(1,0,1)} = x_1 x_3$ | $\Psi_{(1,0,1)} = -[30x_1^2 + 5x_1(5x_2 - 37x_3) - 5x_2^2 + 30x_3^2 + 25x_2x_3 - 28]/28\sqrt{37}$ |
| $\Psi_{(0,2,0)} = \frac{1}{\sqrt{2}}(x_2^2 - 1)$ | $\Psi_{(0,2,0)} = [5x_1^2 + 10(x_3 - 6x_2)x_1 + 180x_2^2 + 5x_3^2 - 60x_2x_3 - 168]/168\sqrt{2}$ |
| $\Psi_{(0,1,1)} = x_2 x_3$ | $\Psi_{(0,1,1)} = [5x_1^2 - 25x_1(x_2 + x_3) - 30x_2^2 - 30x_3^2 + 185x_2x_3 + 28]/28\sqrt{37}$ |
| $\Psi_{(0,0,2)} = \frac{1}{\sqrt{2}}(x_3^2 - 1)$ | $\Psi_{(0,0,2)} = [5x_1^2 + 10x_1(x_2 - 6x_3) + 5x_2^2 + 180x_3^2 - 60x_2x_3 - 168]/168\sqrt{2}$ |

| Case 3: $\rho_{12} = 1/5, \rho_{13} = 2/5, \rho_{23} = 4/5$ | Case 4: $\rho_{12} = -1/5, \rho_{13} = 2/5, \rho_{23} = -4/5$ |
|---|---|
| $\Psi_{(0,0,0)} = 1$ | $\Psi_{(0,0,0)} = 1$ |
| $\Psi_{(1,0,0)} = \frac{\sqrt{5}}{6}(3x_1 + x_2 - 2x_3)$ | $\Psi_{(1,0,0)} = \frac{\sqrt{5}}{6}(3x_1 - x_2 - 2x_3)$ |
| $\Psi_{(0,1,0)} = \frac{1}{2}\sqrt{\frac{5}{21}}(x_1 + 7x_2 - 6x_3)$ | $\Psi_{(0,1,0)} = -\frac{1}{2}\sqrt{\frac{5}{21}}(x_1 - 7x_2 - 6x_3)$ |
| $\Psi_{(0,0,1)} = -\frac{1}{2}\sqrt{\frac{5}{6}}(x_1 + 3x_2 - 4x_3)$ | $\Psi_{(0,0,1)} = -\frac{1}{2}\sqrt{\frac{5}{6}}(x_1 - 3x_2 - 4x_3)$ |
| $\Psi_{(2,0,0)} = [45x_1^2 + 30x_1(x_2 - 2x_3) + 5x_2^2 + 20x_3^2 - 20x_2x_3 - 36]/36\sqrt{2}$ | $\Psi_{(2,0,0)} = [45x_1^2 - 30x_1(x_2 + 2x_3) + 5x_2^2 + 20x_3^2 + 20x_2x_3 - 36]/36\sqrt{2}$ |
| $\Psi_{(1,1,0)} = [15x_1^2 + 10x_1(11x_2 - 10x_3) + 35x_2^2 + 60x_3^2 - 100x_2x_3 - 12]/12\sqrt{22}$ | $\Psi_{(1,1,0)} = -[15x_1^2 - 10x_1(11x_2 + 10x_3) + 35x_2^2 + 60x_3^2 + 100x_2x_3 - 12]/12\sqrt{22}$ |
| $\Psi_{(1,0,1)} = -[15x_1^2 + 10x_1(5x_2 - 7x_3) + 15x_2^2 + 40x_3^2 - 50x_2x_3 - 12]/12\sqrt{7}$ | $\Psi_{(1,0,1)} = -[15x_1^2 - 10x_1(5x_2 + 7x_3) + 15x_2^2 + 40x_3^2 + 50x_2x_3 - 12]/12\sqrt{7}$ |
| $\Psi_{(0,2,0)} = [5x_1^2 + 10x_1(7x_2 - 6x_3) + 245x_2^2 + 180x_3^2 - 420x_2x_3 - 84]/84\sqrt{2}$ | $\Psi_{(0,2,0)} = [5x_1^2 - 10x_1(7x_2 + 6x_3) + 245x_2^2 + 180x_3^2 + 420x_2x_3 - 84]/84\sqrt{2}$ |
| $\Psi_{(0,1,1)} = -[5x_1^2 + 50x_1(x_2 - x_3) + 105x_2^2 + 120x_3^2 - 230x_2x_3 - 36]/12\sqrt{23}$ | $\Psi_{(0,1,1)} = [5x_1^2 - 50x_1(x_2 + x_3) + 105x_2^2 + 120x_3^2 + 230x_2x_3 - 36]/12\sqrt{23}$ |
| $\Psi_{(0,0,2)} = [5x_1^2 + 10x_1(3x_2 - 4x_3) + 45x_2^2 + 80x_3^2 - 120x_2x_3 - 24]/24\sqrt{2}$ | $\Psi_{(0,0,2)} = [5x_1^2 - 10x_1(3x_2 + 4x_3) + 45x_2^2 + 80x_3^2 + 120x_2x_3 - 24]/24\sqrt{2}$ |

Finally, the means and variances of $y_2(\mathbf{X})$ for all four cases, calculated using (35) and (36), are displayed in Table 3. They match the corresponding statistics of $y(\mathbf{X})$ due to the polynomial exactness of the generalized PCE.

## 6.2. Example 2

Consider a stochastic ordinary differential equation (ODE) [13]

$$\frac{dy(t; \mathbf{X})}{dt} = -(1 + X_1)\left[y(t; X_1, X_2) - (1 + X_2)\right], \ 0 \leq t \leq 1,$$

with a deterministic initial condition $y(0; \mathbf{X}) = 0$, where $t$ is an independent variable and $\mathbf{X} = (X_1, X_2)^T$ is a bivariate Gaussian input random vector. The random input has mean $\boldsymbol{\mu}_{\mathbf{X}} = \mathbb{E}[\mathbf{X}] = \mathbf{0} \in \mathbb{R}^2$, positive-definite covariance matrix

$$\mathbf{\Sigma}_{\mathbf{X}} = \mathbb{E}\left[\mathbf{X}\mathbf{X}^T\right] = \begin{bmatrix} \sigma_1^2 & \rho\sigma_1\sigma_2 \\ \rho\sigma_1\sigma_2 & \sigma_2^2 \end{bmatrix} \in \mathbb{S}_+^2,$$

Table 2: Zeroth-, first-, and second-order generalized PCE coefficients in Example 1.

| $C_{(j_1,j_2,j_3)}$ | Case 1 | Case 2 | Case 3 | Case 4 |
|---|---|---|---|---|
| $C_{(0,0,0)}$ | 12 | $\dfrac{63}{5}$ | $\dfrac{67}{5}$ | $\dfrac{57}{5}$ |
| $C_{(1,0,0)}$ | 4 | $2\sqrt{\dfrac{42}{5}}$ | $\dfrac{16}{\sqrt{5}}$ | $\dfrac{12}{\sqrt{5}}$ |
| $C_{(0,1,0)}$ | 4 | $2\sqrt{\dfrac{42}{5}}$ | $4\sqrt{\dfrac{35}{3}}$ | 0 |
| $C_{(0,0,1)}$ | 4 | $2\sqrt{\dfrac{42}{5}}$ | $44\sqrt{\dfrac{2}{15}}$ | $4\sqrt{\dfrac{6}{5}}$ |
| $C_{(2,0,0)}$ | 0 | $\dfrac{33}{35\sqrt{2}}$ | $\dfrac{17}{10\sqrt{2}}$ | $\dfrac{3}{10\sqrt{2}}$ |
| $C_{(1,1,0)}$ | 1 | $\dfrac{19\sqrt{37}}{70}$ | $\dfrac{31}{15}\sqrt{\dfrac{11}{2}}$ | $\dfrac{3}{5}\sqrt{\dfrac{11}{2}}$ |
| $C_{(1,0,1)}$ | 1 | $\dfrac{19\sqrt{37}}{70}$ | $\dfrac{32\sqrt{7}}{15}$ | $-\dfrac{\sqrt{7}}{5}$ |
| $C_{(0,2,0)}$ | 0 | $\dfrac{33}{35\sqrt{2}}$ | $\dfrac{203}{30\sqrt{2}}$ | $-\dfrac{49}{10\sqrt{2}}$ |
| $C_{(0,1,1)}$ | 1 | $\dfrac{19\sqrt{37}}{70}$ | $\dfrac{34\sqrt{23}}{15}$ | $\dfrac{7\sqrt{23}}{5}$ |
| $C_{(0,0,2)}$ | 0 | $\dfrac{33}{35\sqrt{2}}$ | $\dfrac{76\sqrt{2}}{15}$ | $-\dfrac{12\sqrt{2}}{5}$ |

Case 1: $\rho_{12} = \rho_{13} = \rho_{23} = 0$,
Case 2: $\rho_{12} = \rho_{13} = \rho_{23} = 1/5$,
Case 3: $\rho_{12} = 1/5, \rho_{13} = 2/5, \rho_{23} = 4/5$,
Case 4: $\rho_{12} = -1/5, \rho_{13} = 2/5, \rho_{23} = -4/5$.

Table 3: Second-moment properties of $y_2(\mathbf{X})$ in Example 1.

| Case | Mean | Variance |
|---|---|---|
| Case 1 | 12 | 51 |
| Case 2 | $\dfrac{63}{5}$ | $\dfrac{1794}{25}$ |
| Case 3 | $\dfrac{67}{5}$ | $\dfrac{2514}{25}$ |
| Case 4 | $\dfrac{57}{5}$ | $\dfrac{774}{25}$ |

Case 1: $\rho_{12} = \rho_{13} = \rho_{23} = 0$,
Case 2: $\rho_{12} = \rho_{13} = \rho_{23} = 1/5$,
Case 3: $\rho_{12} = 1/5, \rho_{13} = 2/5, \rho_{23} = 4/5$,
Case 4: $\rho_{12} = -1/5, \rho_{13} = 2/5, \rho_{23} = -4/5$.

comprising variances $\sigma_1^2 = \sigma_2^2 = 1/4$ of $X_1$ and $X_2$ and correlation coefficient $-1 < \rho < 1$ between $X_1$ and $X_2$, and a joint density function described by (1) for $N = 2$. The objective of this example is to assess the approximation quality of the truncated generalized PCE in terms of the second-moment statistics of the solution of the ODE.

A direct integration of the stochastic ODE leads to the exact solution: $y(t; \mathbf{X}) = (1 + X_2)[1 - \exp\{-(1 + X_1)t\}]$. As a result, the first two raw moments $\mathbb{E}[y(t; \mathbf{X})]$ and $\mathbb{E}[y^2(t; \mathbf{X})]$, described in Appendix B, can be obtained exactly. Using (B.1) and (B.2), Figure 1 illustrates the plots of the mean $\mathbb{E}[y(t; \mathbf{X})]$ and variance $\mathbb{E}[y^2(t; \mathbf{X})] - (\mathbb{E}[y(t; \mathbf{X})])^2$ of $y(t; \mathbf{X})$ as a function of $t$ for five values of the correlation coefficient: $\rho = -9/10, -1/2, 0, 1/2, 9/10$. Both statistics grow with $t$ regardless of the correlation coefficient as expected. When the correlation coefficient increases, there is a slight uptick in the mean, but the variance rises sharply. Therefore, the second-moment statistics strongly depend on the correlation properties of random input.

Figure 2 depicts nine plots of three second-order standardized multivariate Hermite orthogonal polynomials $\Psi_{(2,0)}(x_1, x_2)$, $\Psi_{(1,1)}(x_1, x_2)$, and $\Psi_{(0,2)}(x_1, x_2)$ for three distinct values of the correlation coefficient: $\rho = -1/2, 0, 1/2$. The polynomials obtained for dependent ($\rho = -1/2, 1/2$) variables are very different than those derived for indepen-dent ($\rho = 0$) variables. Similar plots can be generated for other orders, but they are excluded for brevity.

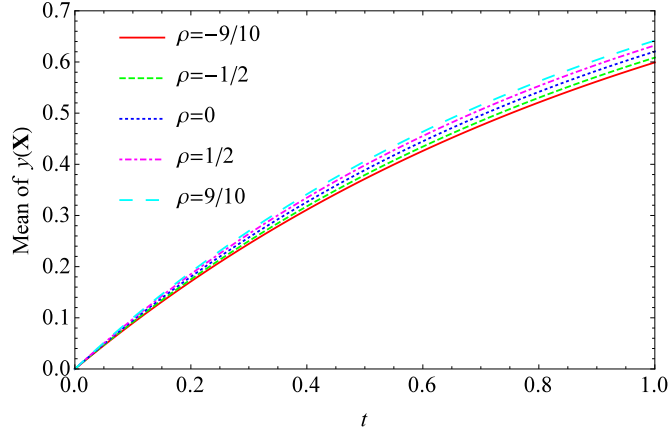

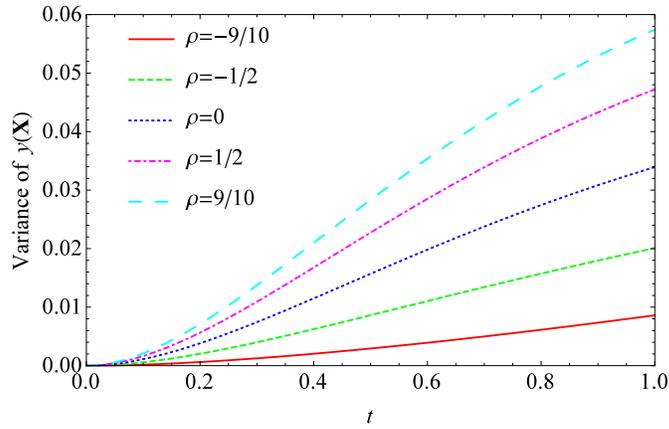

Figure 1: Second-moment statistics of $y(t; \mathbf{X})$ in Example 2; (a) mean; (b) variance.

Since $y(t; \mathbf{X})$ is a non-polynomial function, a convergence analysis with respect to $m$ – the order of the generalized PCE approximation – is essential. Employing $m = 1, 2, 3, 4, 5, 6$ in Algorithm 1, six PCE approximations of $y(t; \mathbf{X})$ and their second-moment statistics were constructed or calculated. Define at $t = 1$ an $L^1$ error

$$e_m := \frac{|\text{var}[y(1; \mathbf{X})] - \text{var}[y_m(1; \mathbf{X})]|}{\text{var}[y(1; \mathbf{X})]} \tag{38}$$

in the variance, committed by an $m$th-order generalized PCE approximation $y_m(1; \mathbf{X})$ of $y(1; \mathbf{X})$, where $\text{var}[y(1; \mathbf{X})]$ and $\text{var}[y_m(1; \mathbf{X})]$ are exact and approximate variances, respectively. The exact variance was obtained from (B.1) and (B.2), whereas the approximate variance, given $m$, was calculated following Algorithm 1. All expectations involved in $\mathbf{A}_l$ and $\mathbf{b}_l$, $0 \le l \le m$, were obtained exactly either by analytical formulae or analytical integrations as in Example 1. Therefore, the variances from the PCE approximations and resultant errors, listed specifically for $\rho = 1/2$ in Appendix B, were determined exactly.

Figure 3 presents five plots describing how the error $e_m$, calculated for each of the five correlation coefficients, decays with respect to $m$. The attenuation rates for all five correlation coefficients are very similar, although the errors for negative correlations are larger than those for non-negative correlations. Nonetheless, nearly exponential conver-gence is achieved by the generalized PCE approximations, preserving the exponential convergence of the classical PCE approximations.

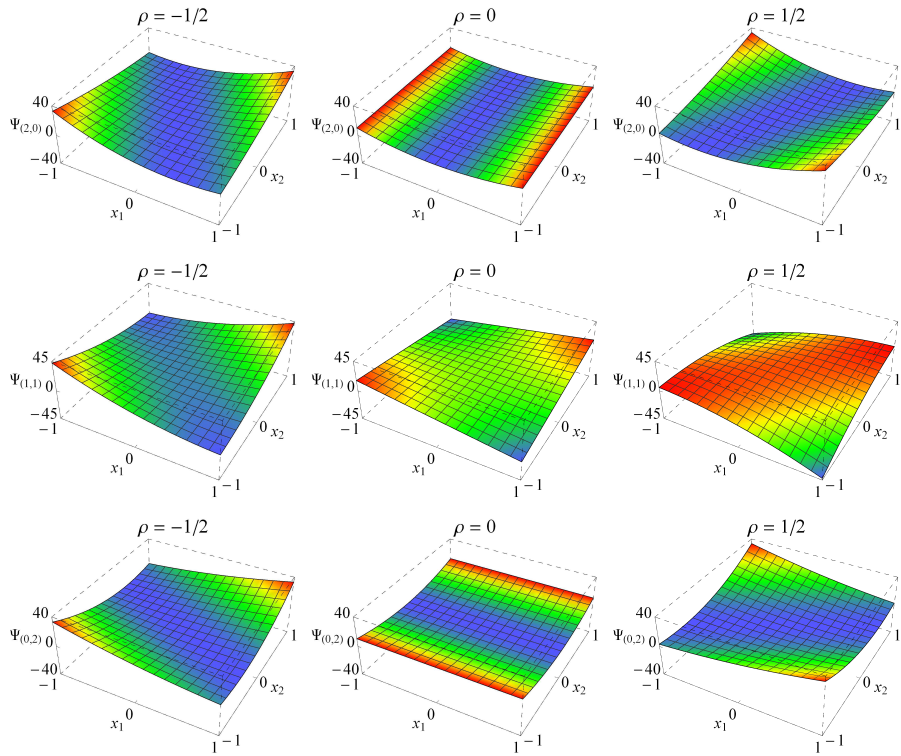

Figure 2: A family of bivariate ($N = 2$), second-order ($m = 2$) standardized Hermite orthogonal polynomials in Example 2 for $\rho = -1/2$ (left), $\rho = 0$ (middle), and $\rho = 1/2$ (right).

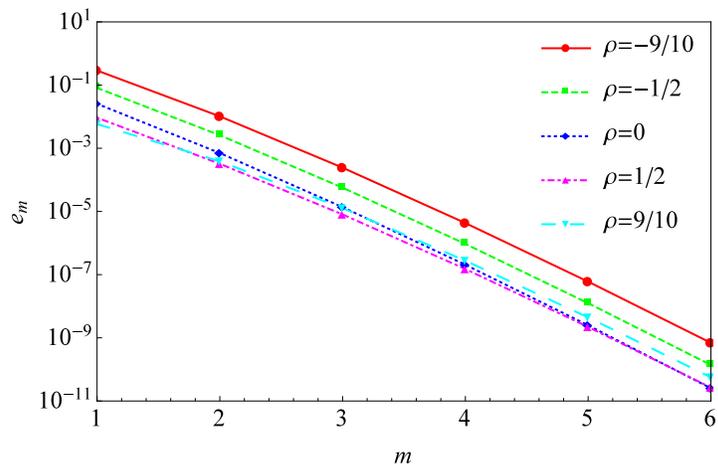

Figure 3: Decay of $L^1$ error in the variance of $y_m(1; \mathbf{X})$ in Example 2 with respect to $m$.

## 6.3. Example 3

The final example entails random eigenvalue analysis of an undamped cantilever plate, shown in Figure 4(a), often performed in structural dynamics. The plate has the following deterministic geometric and material properties: length $L = 2$ in (50.8 mm), width $W = 1$ in (25.4 mm), Young's modulus $E = 30 \times 10^6$ psi (206.8 GPa), Poisson's ratio $\nu = 0.3$, and mass density $\rho = 7.324 \times 10^{-4}$ lb-s$^2$/in$^4$ (7827 kg/mm$^3$). The randomness in eigenvalues arises due to random thickness $t(\xi)$, which is spatially varying in the longitudinal direction $\xi$ only. The thickness is represented by a homogeneous, lognormal random field $t(\xi) = c \exp[\alpha(\xi)]$ with mean $\mu_t = 0.01$ in (0.254 mm), variance $\sigma_t^2 = v_t^2 \mu_t^2$, and coefficient of variation $v_t = 0.2$, where $c = \mu_t / \sqrt{1 + v_t^2}$ and $\alpha(\xi)$ is a zero-mean, homogeneous, Gaussian random field with variance $\sigma_\alpha^2 = \ln(1 + v_t^2)$ and covariance function $\Gamma_\alpha(\tau) = \mathbb{E}[\alpha(\xi)\alpha(\xi + \tau)] = \sigma_\alpha^2 \exp[-|\tau|/(0.2L)]$. Two numerical grids were employed: (1) a $10 \times 20$ finite-element grid of the plate, consisting of 200 eight-noded, second-order shell elements and 661 nodes, as shown in Figure 4(b); and (2) an 11-point random-field grid of the plate, parameterizing the random field $\alpha(\xi)$ into a zero-mean, 11-dimensional, dependent Gaussian random vector $\mathbf{X} = (\alpha_1, \ldots, \alpha_{11})^T$ with covariance matrix $\mathbf{\Sigma_X} = [\Gamma_\alpha(\xi_i - \xi_j)]$, $i, j = 1, \ldots, 11$, where $\xi_i$ is the coordinate of the column of nodes after traversing $2(i-1)$ columns of finite elements from the left, as shown in Figure 4(c). The thickness is linearly interpolated between two consecutive nodes of the random-field grid. The finite-element grid was used for domain discretization, generating the random mass matrix $\mathbf{M(X)}$ and random stiffness matrix $\mathbf{K(X)}$ of the cantilever plate. The random eigenvalue problem calls for solving the matrix characteristic equation: $\det[\mathbf{K(X)} - \mathbf{\Lambda(X)M(X)}] = 0$, where $\mathbf{\Lambda(X)}$ is a random eigenvalue of interest with its square-root representing the corresponding natural frequency. A Lanczos algorithm [31] was used to calculate the eigenvalue.

Using Algorithm 1, the first- and second-order generalized PCE approximations were employed to estimate various probabilistic characteristics of the first four eigenvalues of the plate. The expectations involved in $\mathbf{A}_l$, $l = 0, 1, 2$, were exactly determined from the analytical formulae described by (19) and (20) as before. However, unlike the two former examples, the expectations contained in $\mathbf{b}_l$, $l = 0, 1, 2$, which require 11-dimensional integrations, cannot be determined exactly. Instead, a QMCS was used to estimate the integrals by three steps: (1) select a QMCS sample size $L_{QMCS} \in \mathbb{N}$ and generate a low-discrepancy point set $\mathcal{P}_{L_{QMCS}} := \{\mathbf{u}^{(k)} \in [0, 1]^{11}, k = 1, \ldots, L_{QMCS}\}$; (2) map each sample from $\mathcal{P}_{L_{QMCS}}$ to the sample $\mathbf{x}^{(k)} \in \mathbb{R}^{11}$, following the Gaussian probability measure of $\mathbf{X}$; and (3) approximate the expectation $\mathbb{E}[y(\mathbf{X})\Psi_\mathbf{j}(\mathbf{X}; \mathbf{\Sigma_X})]$ by $\sum_{k=1}^{L_{QMCS}} y(\mathbf{x}^{(k)})\Psi_\mathbf{j}(\mathbf{x}^{(k)}; \mathbf{\Sigma_X})/L_{QMCS}$. The computational cost is proportional to $L_{QMCS}$, as all sample calculations require the same effort. The Sobol sequence [32] was used for the low-discrepancy point set with three distinct values of $L_{QMCS} = 1000, 2000, 3000$.

Table 4: Second-moment properties of first four eigenvalues of the cantilever plate in Example 3.

| Eigenvalue | 1st-order gen. PCE ($L_{QMCS} = 3000$) | | 2nd-order gen. PCE ($L_{QMCS} = 3000$) | | Crude MCS ($L_{MCS} = 10,000$) | |
|---|---|---|---|---|---|---|
| | Mean | St. dev. | Mean | St. dev. | Mean | St. dev. |
| $\Lambda_1$, (rad/ms)$^2$ | 0.275088 | 0.0869715 | 0.275088 | 0.0896882 | 0.274852 | 0.0888108 |
| $\Lambda_2$, (rad/ms)$^2$ | 5.10714 | 1.18771 | 5.10714 | 1.21458 | 5.10376 | 1.20242 |
| $\Lambda_3$, (rad/ms)$^2$ | 10.6004 | 2.0924 | 10.6004 | 2.14212 | 10.5987 | 2.14294 |
| $\Lambda_4$, (rad/ms)$^2$ | 54.5265 | 9.85134 | 54.5265 | 10.1103 | 54.5506 | 10.0225 |

Table 4 presents the means and standard deviations of the first four eigenvalues, $\Lambda_i$, $i = 1, \ldots, 4$, of the plate by three different methods: the two generalized PCE approximations and crude MCS. The expansion coefficients of the PCE are based on the QMCS sample size $L_{QMCS} = 3000$. In all three methods, the solution of the matrix characteristic equation for a given input is equivalent to performing a finite-element analysis. Therefore, computational efficiency, even for this simple plate, is a practical requirement in solving random eigenvalue problems. Due to the expense of finite-element analysis, crude MCS was conducted for a sample size $L_{MCS} = 10,000$, which should be adequate for providing benchmark solutions of the second-moment characteristics. The agreement between the means and standard deviations by both PCE approximations and crude MCS in Table 4 is good. However, the second-order approximation is relatively more accurate than the first-order approximation in estimating standard deviations, as expected.

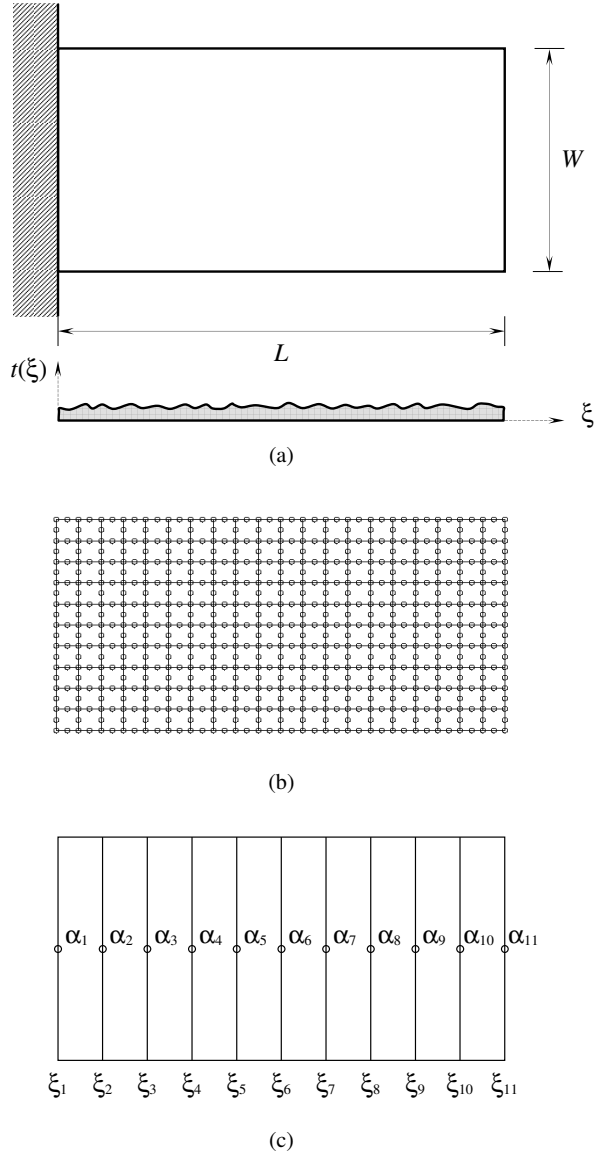

Figure 4: A cantilever plate; (a) geometry; (b) finite-element grid; (c) random-field grid.

Figures 5 and 6 illustrate the marginal probability density functions of the four eigenvalues by the two generalized PCE approximations and crude MCS. Due to the computational expense inherent to finite-element analysis, the same 10,000 samples generated for verifying the statistics in Table 4 were utilized to develop the histograms of crude MCS in Figures 5 and 6. However, since the PCE approximations yield explicit eigenvalue approximations in terms of multivariate polynomials, a relatively large sample size, 100,000 in this particular example, was selected to sample (33) for estimating the respective densities by histograms as well. Moreover, for each eigenvalue and order, three PCE-based densities, obtained when estimating the expansion coefficients with $L_{QMCS}$ = 1000, 2000, and 3000, were generated to monitor convergence. The respective densities estimated by the second-order PCE approximations and crude MCS match well over the entire support for all four eigenvalues, especially when the QMCS sample size is relatively large. In contrast, the first-order PCE approximations produce satisfactory density estimates only around the means; there are discrepancies in the tail regions of the densities even when $L_{QMCS} = 3000$. This suggests that a satisfactory second-moment analysis by the first-order PCE approximation may not translate to accurate calculation of the probability density function. This known problem for the classical PCE persists for the generalized PCE.

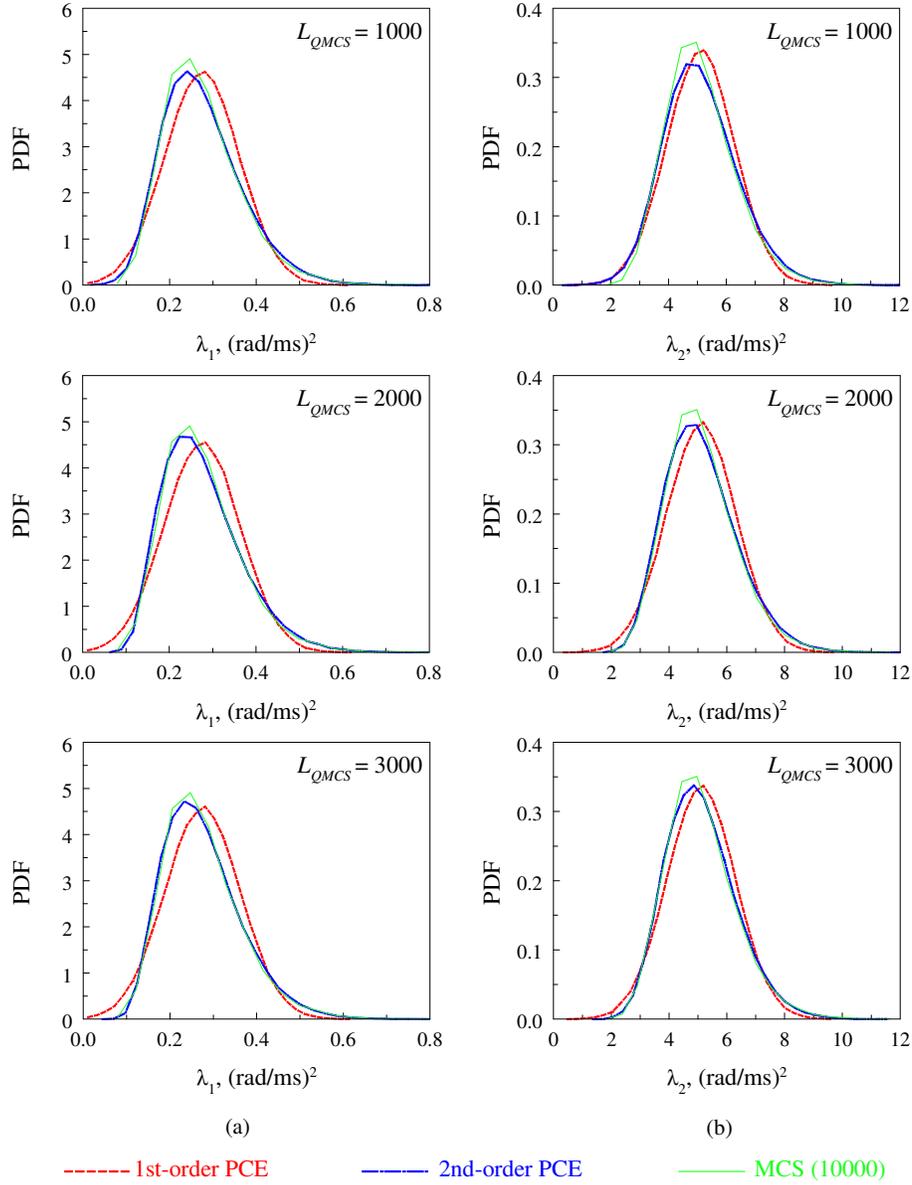

Figure 5: Marginal probability density functions (PDFs) of random eigenvalues of the cantilever plate in Example 3 by generalized PCE approximations and crude MCS; (a) first eigenvalue; (2) second eigenvalue.

## 7. Conclusion

A new generalized PCE of a square-integrable random variable, comprising multivariate Hermite polynomials in dependent Gaussian random variables, is presented. Derived analytically, the second-moment properties of multivariate Hermite polynomials reveal a weakly orthogonal system with respect to an inner product comprising a general Gaussian probability measure. When the Gaussian variables are statistically independent, the multivariate Hermite polynomials elevate to a strongly orthogonal system, leading to the classical PCE. Nonetheless, when the Gaussian variables are statistically dependent, the exponential integrability of norm still allows the Hermite polynomials to constitute a complete set and hence a basis in a Hilbert space. The completeness is vitally important for the convergence of the generalized PCE to the correct limit. The optimality of the generalized PCE and the approximation quality due to truncation have been discussed. New analytical formulae are proposed to calculate the mean and variance of a generalized PCE approximation of a general output variable in terms of the expansion coefficients and statistical properties of Hermite polynomials. However, unlike in the classical PCE, calculating the coefficients of the generalized PCE requires solving a coupled system of linear equations. Moreover, the variance formula of the generalized PCE contains additional terms – a consequence of statistical dependence among Gaussian input variables – that are not present in that of the classical PCE. The additional terms vanish as they should when the Gaussian variables are statistically independent, regressing the generalized PCE to the classical PCE. A possible extension of the generalized PCE for non-Gaussian variables has been discussed. Numerical examples developed from an elementary function, a stochastic ODE, and a random eigenvalue analysis illustrate the construction and use of a generalized PCE approximation in estimating the statistical properties of output variables.

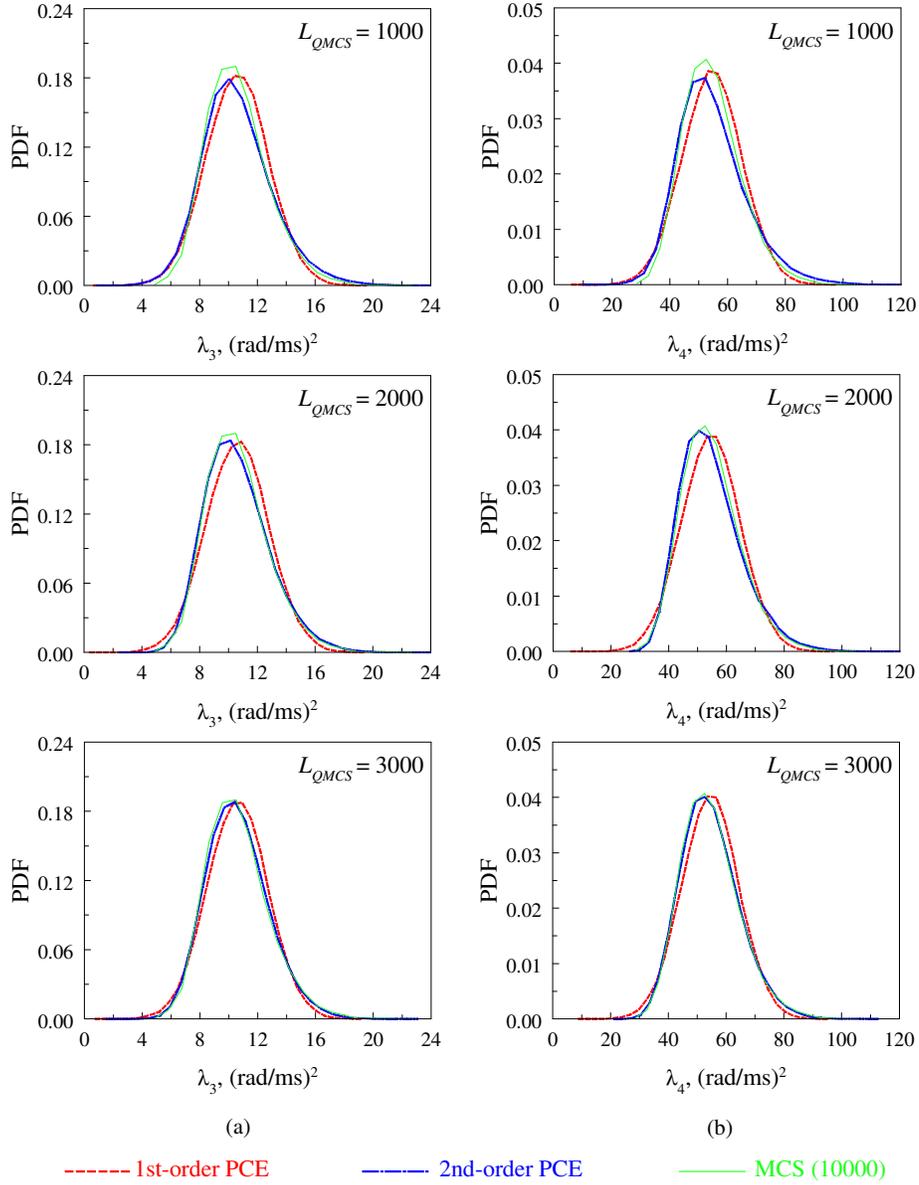

Figure 6: Marginal probability density functions (PDFs) of random eigenvalues of the cantilever plate in Example 3 by generalized PCE approximations and crude MCS; (a) third eigenvalue; (2) fourth eigenvalue.

## Appendix A. Matrix form for calculating expansion coefficients

The linear system (37) for $|\mathbf{j}| = l$ involves $K_{N,l}$ number of orthornormal polynomials and expansion coefficients, which must be ordered with a single index, say, $p$. For multivariate polynomials, there are many options. One option employed in this work is the graded lexicographic order.

**Definition 18.** *A general monomial order is denoted by the symbol $>$. For $\mathbf{j}, \mathbf{k} \in \mathbb{N}_0^N$, the graded lexicographic order, denoted by $>_{grlex}$, is such that $\mathbf{j} >_{grlex} \mathbf{k}$ if and only if $|\mathbf{j}| \geq |\mathbf{k}|$ and the leftmost nonzero entry of $\mathbf{j} - \mathbf{k}$ is positive.*

Using the graded lexicographic order from Definition 18, the multi-indices $\mathbf{j}$, $|\mathbf{j}| = l$, can now be arranged in an ascending order following a single index $p$, which runs from 1 to $K_{N,l}$. Table A.5 illustrates the graded lexicographic order for a three-dimensional case ($N = 3$) and three subcases: $l = 0$, $l = 1$, and $l = 2$.

Table A.5: Graded lexicographic order of multi-index $\mathbf{j}$ for $N = 3$. and $l = 0, 1, 2$

| $|\mathbf{j}| = l$ | $K_{3,l}$ | Multi-index $\mathbf{j}$ | Single index $p$ |
|---|---|---|---|
| 0 | 1 | (0,0,0) | 1 |
| 1 | 3 | (1,0,0) | 1 |
| | | (0,1,0) | 2 |
| | | (0,0,1) | 3 |
| 2 | 6 | (2,0,0) | 1 |
| | | (1,1,0) | 2 |
| | | (1,0,1) | 3 |
| | | (0,2,0) | 4 |
| | | (0,1,1) | 5 |
| | | (0,0,2) | 6 |

The matrix form of (37) requires construction of the following: (1) a $K_{N,l} \times K_{N,l}$ matrix $\mathbf{A}_l \in \mathbb{S}_+^{K_{N,l}}$, comprising the expectations $\mathbb{E}[\Psi_\mathbf{j}(\mathbf{X}; \Sigma_\mathbf{X})\Psi_\mathbf{k}(\mathbf{X}; \Sigma_\mathbf{X})]$, $|\mathbf{j}| = |\mathbf{k}| = l$; (2) a $K_{N,l}$-dimensional vector $\mathbf{b}_l \in \mathbb{R}^{K_{N,l}}$, consisting of the expectations $\mathbb{E}[y(\mathbf{X})\Psi_\mathbf{j}(\mathbf{X}; \Sigma_\mathbf{X})]$, $|\mathbf{j}| = l$; and (3) a $K_{N,l}$-dimensional vector $\mathbf{c}_l \in \mathbb{R}^{K_{N,l}}$, collecting the expansion coefficients $C_\mathbf{j}$, $|\mathbf{j}| = l$. The elements of the system matrix and vectors are arranged according to the graded lexicographic order described earlier. For example, when $N = 3$ and $l = 2$, the size of the linear system is $K_{3,2} = (4!)/(2!2!) = 6$, yielding

$$\mathbf{A}_2 = \begin{bmatrix} \mathbb{E}[\Psi_{(2,0,0)}^2] & \mathbb{E}[\Psi_{(2,0,0)}\Psi_{(1,1,0)}] & \mathbb{E}[\Psi_{(2,0,0)}\Psi_{(1,0,1)}] & \mathbb{E}[\Psi_{(2,0,0)}\Psi_{(0,2,0)}] & \mathbb{E}[\Psi_{(2,0,0)}\Psi_{(0,1,1)}] & \mathbb{E}[\Psi_{(2,0,0)}\Psi_{(0,0,2)}] \\ & \mathbb{E}[\Psi_{(1,1,0)}^2] & \mathbb{E}[\Psi_{(1,1,0)}\Psi_{(1,0,1)}] & \mathbb{E}[\Psi_{(1,1,0)}\Psi_{(0,2,0)}] & \mathbb{E}[\Psi_{(1,1,0)}\Psi_{(0,1,1)}] & \mathbb{E}[\Psi_{(1,1,0)}\Psi_{(0,0,2)}] \\ & & \mathbb{E}[\Psi_{(1,0,1)}^2] & \mathbb{E}[\Psi_{(1,0,1)}\Psi_{(0,2,0)}] & \mathbb{E}[\Psi_{(1,0,1)}\Psi_{(0,1,1)}] & \mathbb{E}[\Psi_{(1,0,1)}\Psi_{(0,0,2)}] \\ & & & \mathbb{E}[\Psi_{(0,2,0)}^2] & \mathbb{E}[\Psi_{(0,2,0)}\Psi_{(0,1,1)}] & \mathbb{E}[\Psi_{(0,2,0)}\Psi_{(0,0,2)}] \\ & & & & \mathbb{E}[\Psi_{(0,1,1)}^2] & \mathbb{E}[\Psi_{(0,1,1)}\Psi_{(0,0,2)}] \\ (\text{sym.}) & & & & & \mathbb{E}[\Psi_{(0,0,2)}^2] \end{bmatrix},$$

$$\mathbf{b}_2 = \begin{pmatrix} \mathbb{E}[y\Psi_{(2,0,0)}] \\ \mathbb{E}[y\Psi_{(1,1,0)}] \\ \mathbb{E}[y\Psi_{(1,0,1)}] \\ \mathbb{E}[y\Psi_{(0,2,0)}] \\ \mathbb{E}[y\Psi_{(0,1,1)}] \\ \mathbb{E}[y\Psi_{(0,0,2)}] \end{pmatrix}, \quad \mathbf{c}_2 = \begin{pmatrix} C_{(2,0,0)} \\ C_{(1,1,0)} \\ C_{(1,0,1)} \\ C_{(0,2,0)} \\ C_{(0,1,1)} \\ C_{(0,0,2)} \end{pmatrix}.$$

From Corollary 12, $\mathbf{A}_l \in \mathbb{R}^{K_{N,l} \times K_{N,l}}$ is a symmetric, positive-definite matrix, and hence invertible. The solution of $\mathbf{A}_l \mathbf{c}_l = \mathbf{b}_l$ produces the expansion coefficients.

### Appendix B. Exact second-moment properties of $y(t; \mathbf{X})$ and approximation errors in variance

Applying the expectation operators on the solution of the stochastic ODE from Example 2 and its square, the first two raw moments of $y(t; \mathbf{X})$, valid for $t \in [0, 1]$ and $\rho \in (-1, 1)$, respectively, are

$$\mathbb{E}[y(t; \mathbf{X})] = 1 + \left(\frac{\rho t}{16} - 1\right) \exp\left(\frac{t^2}{32} - t\right) \tag{B.1}$$

and

$$\begin{aligned} \mathbb{E}\left[y^2(t; \mathbf{X})\right] &= \frac{1}{128} \exp(-2t) \left[ 136 \exp(2t) - \exp\left(t + \frac{t^2}{32}\right) \{272 + \rho t(\rho t - 32)\} \right. \\ &\quad \left. + 2 \exp\left(\frac{t^2}{8}\right) \{68 + \rho t(\rho t - 16)\} \right]. \end{aligned} \tag{B.2}$$

From (38), six PCE committed $L^1$ approximation errors in variance for $\rho = 1/2$ are

$$e_1 = \frac{15424 \sqrt[16]{e} - 16369}{16\left(-961 + 964 \sqrt[16]{e} - 66 e^{31/32} + 64 e^{31/16}\right)} \approx 9.26928 \times 10^{-3},$$

$$e_2 = \frac{493568 \sqrt[16]{e} - 525345}{512\left(-961 + 964 \sqrt[16]{e} - 66 e^{31/32} + 64 e^{31/16}\right)} \approx 3.22487 \times 10^{-4},$$

$$e_3 = \frac{23691264 \sqrt[16]{e} - 25219153}{24576\left(-961 + 964 \sqrt[16]{e} - 66 e^{31/32} + 64 e^{31/16}\right)} \approx 8.03445 \times 10^{-6},$$

$$e_4 = \frac{1516240896 \sqrt[16]{e} - 1614029953}{1572864\left(-961 + 964 \sqrt[16]{e} - 66 e^{31/32} + 64 e^{31/16}\right)} \approx 1.50027 \times 10^{-7},$$

$$e_5 = \frac{40433090560 \sqrt[16]{e} - 43040800827}{41943040\left(-961 + 964 \sqrt[16]{e} - 66 e^{31/32} + 64 e^{31/16}\right)} \approx 2.20588 \times 10^{-9},$$

$$e_6 = \frac{11644730081280 \sqrt[16]{e} - 12395750647009}{12079595520\left(-961 + 964 \sqrt[16]{e} - 66 e^{31/32} + 64 e^{31/16}\right)} \approx 2.65667 \times 10^{-11}.$$

Similar results, also generated for other values of $\rho$, are not reported here for brevity.

### Acknowledgments

The author sincerely thanks the division editor and reviewers for providing many useful comments.